\documentclass[a4paper,11pt,reqno]{amsart}

\usepackage[T1]{fontenc}
\usepackage[french,english]{babel}

\usepackage[headings]{fullpage}

\usepackage{dsfont,amssymb}
\usepackage{tikz}
\usetikzlibrary{cd}
\usepackage[all]{xy}

\usepackage{hyperref}

\usepackage[nobysame,alphabetic,initials,msc-links]{amsrefs}
\DefineSimpleKey{bib}{how}

\renewcommand{\eprint}[1]{#1}
\BibSpec{misc}{%
  +{}{\PrintAuthors}  {author}
  +{,}{ \textit}      {title}
  +{,}{ }             {how}
  +{}{ \parenthesize} {date}
  +{,} { available at \eprint}        {eprint}
  +{,}{ available at \url}{url}
  +{,}{ }             {note}
  +{.}{}              {transition}
}

\mathchardef\mhyph="2D

\numberwithin{equation}{section}

\newtheorem{theorem}{Theorem}[section]
\newtheorem{corollary}[theorem]{Corollary}
\newtheorem{lemma}[theorem]{Lemma}
\newtheorem{proposition}[theorem]{Proposition}
\theoremstyle{remark}
\newtheorem{remark}[theorem]{Remark}
\newtheorem{example}[theorem]{Example}
\theoremstyle{definition}
\newtheorem{definition}[theorem]{Definition}

\newcommand\bp{\begin{proof}}
\newcommand\ep{\end{proof}}

\usepackage{wasysym}
\newcommand{\circt}%
{\mathbin{%
\mathchoice
{\ooalign{$\ocircle$\cr\hidewidth\raise-.15ex\hbox{$\scriptstyle\top\mkern2.05mu$}\cr}}% Woronowicz style tensor product, USUAL SIZE
{\ooalign{$\ocircle$\cr\hidewidth\raise-.15ex\hbox{$\scriptstyle\top\mkern2.05mu$}\cr}}% Woronowicz style tensor product, USUAL SIZE
{\ooalign{$\scriptstyle\ocircle$\cr\hidewidth\raise-.12ex\hbox{$\scriptscriptstyle\top\mkern1mu$}\cr}}% Woronowicz style tensor product, SCRIPT SIZE
{\ooalign{$\scriptstyle\ocircle$\cr\hidewidth\raise-.12ex\hbox{$\scriptscriptstyle\top\mkern1mu$}\cr}}% Woronowicz style tensor product, SCRIPT SIZE
}}

\DeclareMathOperator{\Ad}{Ad}
\DeclareMathOperator{\QAut}{QAut}
\DeclareMathOperator{\End}{End}
\DeclareMathOperator{\Hilbf}{\operatorname{Hilb}_{\mathrm f}}
\DeclareMathOperator{\Mat}{Mat}
\DeclareMathOperator{\Mor}{Mor}
\DeclareMathOperator{\Nat}{Nat}
\DeclareMathOperator{\Irr}{Irr}
\DeclareMathOperator{\Rep}{Rep}
\DeclareMathOperator{\sign}{sign}
\DeclareMathOperator{\Tr}{Tr}
\newcommand\tr{\operatorname{tr}}
\newcommand\Dhat{\hat\Delta}

\newcommand\algg{\mathrm{alg}}

\newcommand{\C}{{\mathbb C}}
\newcommand{\N}{{\mathbb N}}

\newcommand\PPP{{\mathbb P}}

\newcommand\T{{\mathbb T}}

\newcommand{\A}{{\mathcal A}}
\newcommand{\B}{{\mathcal B}}

\newcommand\CC{{\mathcal C}}
\newcommand\D{\mathcal D}

\newcommand{\NN}{{\mathcal N}}
\newcommand\RR{\mathcal R}
\newcommand\PP{{\mathcal P}}
\newcommand\U{\mathcal U}
\newcommand\un{{\mathds 1}}
\newcommand\eps{\varepsilon}

\newcommand{\SU}{\mathrm{SU}}
\newcommand{\SO}{\mathrm{SO}}
\newcommand{\SL}{\mathrm{SL}}

\newcommand\AUF{U^+_F}
\newcommand\FAUF{{\mathbb F}U_F}

\newcommand\FH{\mathrm{FH}}
\newcommand\cb{\mathrm{cb}}

\begin{document}

\title[Noncommutative boundaries]{Noncommutative Poisson boundaries and Furstenberg--Hamana boundaries of Drinfeld doubles}

\date{May 7, 2021; minor changes December 27, 2021}

\author{Erik Habbestad}
\address{University of Oslo, Mathematics institute}
\email{erikhab@math.uio.no}

\author{Lucas Hataishi}
%\address{Universitetet i Oslo}
\email{lucasyh@math.uio.no}

\author{Sergey Neshveyev}
%\address{Universitetet i Oslo}
\email{sergeyn@math.uio.no}

\thanks{Supported by the NFR project 300837 ``Quantum Symmetry''.}



\maketitle

%\begin{minipage}{\linewidth}
\begin{abstract}
We clarify the relation between noncommutative Poisson boundaries and Fursten\-berg--Hamana boundaries of quantum groups. Specifically, given a compact quantum group~$G$, we show that in many cases where the Poisson boundary of the dual discrete quantum group $\hat G$ has been computed, the underlying topological boundary either coincides with the Furstenberg--Hamana boundary of the Drinfeld double $D(G)$ of $G$ or is a quotient of it. This includes the $q$-deformations of compact Lie groups, free orthogonal and free unitary quantum groups, quantum automorphism groups of finite dimensional C$^*$-algebras. In particular, the boundary of $D(G_q)$ for the $q$-deformation of a compact connected semisimple Lie group $G$ is~$G_q/T$ (for $q\ne1$), in agreement with the classical results of Furstenberg and Moore on the Furstenberg boundary of $G_\C$.

We show also that the construction of the Furstenberg--Hamana boundary of $D(G)$ respects monoidal equivalence and, in fact, can be carried out entirely at the level of the representation category of $G$. This leads to a notion of the Furstenberg--Hamana boundary of a rigid C$^*$-tensor category.
\end{abstract}

\selectlanguage{french}
\begin{abstract}
Nous clarifions la relation entre fronti\`{e}res de Poisson non commutatives et  fronti\`{e}res de Furstenberg--Hamana des groupes quantiques. Plus pr\'ecis\'ement, \'etant donn\'e un groupe quantique compact~$G$, nous montrons que dans les nombreux cas o\`{u} la fronti\`{e}re de Poisson du groupe quantique dual discret  $\hat G$ a \'et\'e calcul\'e,  la fronti\`{e}re topologique sous-jacente co\"{\i}ncide soit avec la fronti\`{e}re de Furstenberg--Hamana du double de Drinfeld $D(G)$ de $G$ soit est un quotient de celle-ci. Cela inclut les $q$-d\'eformations des groupes de Lie compacts, les groupes quantiques libres orthogonaux et unitaires, les groupes d'automorphismes quantiques des C$^*$-alg\`{e}bres de dimension finie. En particulier, la fronti\`{e}re de $D(G_q)$ pour la $q$-d\'eformation d'un groupe de Lie semi-simple connexe compact $G$ est~$G_q/T$ (pour $q\ne1$), en accord avec les r\'esultats classiques de Furstenberg et Moore sur la fronti\`{e}re de Furstenberg de $G_\C$.

Nous montrons aussi que la construction de la fronti\`{e}re de Furstenberg--Hamana de $D(G)$ respecte l'\'equivalence mono\"{\i}dale et, en fait, peut \^{e}tre r\'ealis\'ee enti\`{e}rement au niveau de la cat\'egorie de repr\'esentation de $G$. Cela conduit \`{a} une notion de  fronti\`{e}re de Furstenberg--Hamana d'une C$^*$-cat\'{e}gorie tensorielle rigide.
\end{abstract}
\selectlanguage{english}
%\end{minipage}

\section*{Introduction}

In his celebrated work on the Poisson formula for semisimple Lie groups~\cite{MR146298} Furstenberg attached two boundaries to every locally compact group $G$, which are now called the Poisson and Furstenberg boundaries of $G$. For real semisimple Lie groups with finite centers he showed that the two boundaries coincide and give rise to an integral representation of the bounded harmonic functions on $G$. As has been observed since then both constructions can be phrased in operator algebraic terms, at least when $G$ is discrete, paving the way to their generalizations to the noncommutative setting.

For the Furstenberg boundary $\partial_{\mathrm F} G$, it was noticed  by Hamana~\cite{MR509025} that $C(\partial_{\mathrm F} G)$ coincides with the injective envelope of the $G$-algebra $\C$. In fact, Hamana's motivation was quite different from that of Furstenberg and the connection between the two works was mentioned only in passing. This connection had not attracted any attention until a few years ago, when it was rediscovered by Kalantar and Kennedy in their work on C$^*$-simple groups~\cite{MR3652252}. Since then Hamana's construction of injective envelopes has been used to develop a Furstenberg boundary-type theory in several different contexts - for unitary representations of discrete groups~\cite{BK}, \'etale groupoids~\cite{Bor}, discrete quantum groups~\cite{KKSV}.

On the other hand, a construction of the Poisson boundary in the operator algebraic framework was given by Izumi~\cite{MR1916370}. He initiated the study of Poisson boundaries of discrete quantum groups and computed the boundary of the dual of $\SU_q(2)$. The answer - the standard Podle\'s quantum sphere $S^2_q$ - turned out to be a quantization of the Poisson boundary of the complexification $\SL_2(\C)$ of $\SU(2)$. The same phenomenon was then demonstrated for~$\SU_q(n)$~\cite{MR2200270}, at which point it became clear that this is not a coincidence and requires an explanation~\cite{Nreview}. A computation of the Poisson boundaries of the duals of all $q$-deformed compact semisimple Lie groups was done soon afterwards by Tomatsu~\cite{MR2335776}.

A satisfactory conceptual explanation of the above phenomenon is actually not difficult to find. Given a compact semisimple Lie group $G$, the Poisson boundary of $\hat G_q$ carries an action of the Drinfeld double $D(G_q)$. It can then be shown (see Proposition~\ref{poisson-boundary-D(G)}) that the Poisson boundary of $D(G_q)$ is isomorphic, as a noncommutative $D(G_q)$-space, to that of $\hat G_q$. Since $D(G_q)$ can be viewed as a quantization of $G_\C$, it is therefore not so surprising that the Poisson boundary of $\hat G_q$ is a quantization of that of $G_\C$. In hindsight, even several steps in the two computations are related. For example, an important property both in the classical~\cite{MR146298} and quantum~\cite{MR1753202} (see also~\cite{MR2200270}) cases is that the Poisson boundary is a homogeneous space of $G$ and $G_q$, resp.

It is then natural to ask what the Furstenberg boundary of $D(G_q)$, and possibly of some other quantum groups, is. In this paper we answer this question using already known properties of noncommutative Poisson boundaries. It is curious that this direction - from Poisson to Furstenberg boundaries - is opposite to the one in Furstenberg's work, but we leave a detailed comparison of the two approaches to another occasion.

%\smallskip

In more detail, the contents of the paper is as follows. After a short preliminary Section~\ref{sec:prelim},  in Section~\ref{sec:FH} we define the Furstenberg--Hamana boundary $C(\partial_\FH D(G))$ of the Drinfeld double~$D(G)$ of a compact quantum group~$G$ as the $D(G)$-injective envelope of $\C$. The existence and uniqueness of this object (Theorem~\ref{thm:FH}) are proved similarly to Hamana's work and subsequent papers. It is worth stressing though that we make use of special properties of the Drinfeld doubles and do not attempt to work with arbitrary locally compact quantum groups (cf.~\cite{MR2985658}). Our setting is essentially a $G$-equivariant version of the case of discrete quantum groups studied recently in~\cite{KKSV}.

In Section~\ref{sec:PvsF} we show that basic properties of the Poisson and Furstenberg--Hamana boundaries quickly imply that if the action of $G$ on the Poisson boundary of $\hat G$ is ergodic, then this boundary coincides with the Furstenberg--Hamana boundary of $D(G)$. (To be more precise, a noncommutative Poisson boundary is in general only a measure-theoretic object, but under the ergodicity assumption there is no question what the underlying topological structure is.) The quantum groups $G$ with this property are precisely the ones with countable isomorphism classes of irreducible representations and weakly amenable quantum dimension functions. The corresponding Poisson boundaries have been shown to be universal objects measuring how far these dimension functions are from amenable ones~\cite{MR3679617}. Using this property we define a Poisson-type boundary $\partial_\Pi\hat G$ of $\hat G$ for any compact quantum group $G$ with weakly amenable quantum dimension function (Theorem~\ref{thm:NY-extension}), and then show that we still have $\partial_\FH D(G)=\partial_\Pi\hat G$ (Theorem~\ref{thm:FH-weak-amenable}). We draw some consequences of this equality and illustrate it with several examples. Going beyond the weakly amenable case, in Section~\ref{ssec:free-unitary} we consider a free unitary quantum group~$U^+_F$ and show that the topological boundary of its dual defined in~\cite{MR2660685} is a quotient of~$\partial_\FH D(U^+_F)$ (Theorem~\ref{thm:free-unitary}).

One of the consequences of the construction of the $D(G)$-algebra $C(\partial_\FH D(G))$ is that it is braided-commutative (Corollary~\ref{cor:br-com}), which by results of~\cite{MR3291643} implies that there is a C$^*$-tensor category associated with it. It is natural to expect that this category depends only on the representation category of~$G$, similarly to the cases of noncommutative Poisson and Martin boundaries~\citelist{\cite{MR2664313}\cite{MR3291643}\cite{MR3725499}}. In Section~\ref{sec:cat} we show that this is indeed the case and define Furstenberg--Hamana boundaries of arbitrary rigid C$^*$-tensor categories with simple units. We explain, without going into too many details, how a number of results from Section~\ref{sec:PvsF} generalize to this setting. But the main new technical point of Section~\ref{sec:cat}, which is of independent interest, is a categorical description of $G$- and $D(G)$-equivariant completely positive, completely bounded and completely isometric maps (Propositions~\ref{prop:multipliers} and~\ref{prop:multipliers2}). Particular cases of this description have appeared in~\cite{MR3679617} to define categorical analogues of invariant means and in the work of Popa and Vaes~\cite{MR3406647}, which developed one of the equivalent approaches~\citelist{\cite{MR3406647}\cite{MR3509018}\cite{MR3447719}} to defining a maximal C$^*$-norm on the fusion algebra of a rigid C$^*$-tensor category. Furthermore, categorical analogues of completely positive maps have been already proposed in full generality in~\cite{MR3687214}, so the definitions we give are natural and essentially known. The crux of the matter is to show that they indeed reflect properties of equivariant maps under a Tannaka--Krein duality for quantum group actions.

\section{Preliminaries}\label{sec:prelim}

\subsection{Compact quantum groups and their duals}
\label{sec:comp-quant-groups}

We follow the conventions of~\cite{neshveyev-tuset-book}, but we recall some of the main concepts for the reader's convenience.
Let $G$ be a compact quantum group and $(\C[G],\Delta)$ be the Hopf $*$-algebra of regular functions on~$G$. Throughout the entire paper we will work only with the reduced form of $G$, so $C(G)$ denotes the closure of $\C[G]$ in the GNS-representation defined by the Haar state $h$.

A finite dimensional unitary representation of $G$ is a unitary element $U$ of $B(H_U)\otimes C(G)$ (or equivalently, of $B(H_U)\otimes \C[G]$), where~$H_U$ is a finite dimensional Hilbert space, such that
$$
(\iota\otimes\Delta)(U)=U_{12}U_{13}.
$$
The tensor product of two representations $U$ and $V$ is defined by $U_{13}V_{23}$ and denoted by $U\circt V$ or simply by $U\otimes V$, when there is no danger of confusing it with the tensor product of operators.  The C$^*$-tensor category of finite dimensional unitary representations of~$G$ is denoted by $\Rep G$. The unit $\un$ of $\Rep G$ is the trivial representation $1\in C(G)$. We denote by $\Irr(G)$ the set of equivalence classes of irreducible unitary representations of $G$. For every $s\in\Irr(G)$ we fix a representative $U_s$ and write $H_s$ for $H_{U_s}$.

\smallskip

The dual space $\U(G)=\C[G]^*$ has the structure of a $*$-algebra, defined by duality from the Hopf $*$-algebra $(\C[G],\Delta)$. We also define $\U(G^n)=(\C[G]^{\otimes n})^*$. Then the dual of the product map $\C[G]\otimes\C[G]\to\C[G]$ is a $*$-homomorphism $\Dhat\colon\U(G)\to \U(G\times G)$. Every finite dimensional unitary representation $U$ of $G$ defines a $*$-representation $\pi_U$ of $\U(G)$ on $H_U$ by $\pi_U(\omega)=(\iota\otimes\omega)(U)$. The representations $\pi_s=\pi_{U_s}$, $s\in\Irr(G)$, allow us to identify $\U(G)$ with $\prod_{s\in\Irr(G)}B(H_s)$. Then the dual discrete quantum group $\hat G$, in the von Neumann algebra setting, is described by the von Neumann algebra
$$
\ell^\infty(\hat G)=\ell^\infty\text{-}\bigoplus_{s\in\Irr(G)}B(H_s)\subset\U(G)
$$
with comultiplication $\Dhat|_{\ell^\infty(\hat G)}$. We also define
$$
c_0(\hat G)=c_0\text{-}\bigoplus_{s\in\Irr(G)}B(H_s),\qquad c_c(\hat G)=\bigoplus_{s\in\Irr(G)}B(H_s).
$$

The fundamental unitary of $G$ is defined by
\begin{equation*}
W=\bigoplus_{s\in\Irr(G)}(U_s)_{21}\in M(C(G)\otimes c_0(\hat G)).
\end{equation*}
We have the following identities, which reflect the duality between $G$ and $\hat G$ and are equivalent to the pentagon relation for $W$ in the regular representations of $G$ and $\hat G$:
\begin{equation}\label{eq:pentagon}
(\Delta\otimes\iota)(W)=W_{13}W_{23},\qquad (\iota\otimes\Dhat)(W)=W_{12}W_{13}.
\end{equation}

\smallskip

The Woronowicz character $f_1\in\U(G)$ is denoted by~$\rho$. Then, given a finite dimensional unitary representation $U$ of $G$, the conjugate unitary representation is defined by
$$
\bar U=(j(\rho_U)^{1/2}\otimes1)(j\otimes\iota)(U^*)(j(\rho_U)^{-1/2}\otimes1)\in B(\bar H_U)\otimes\C[G],
$$
where $\rho_U=\pi_U(\rho)$ and $j$ denotes the canonical $*$-anti-isomorphism $B(H_U)\cong B(\bar H_U)$ defined by $j(T)\bar\xi=\overline{T^*\xi}$. We have morphisms $R_U\colon \un\to \bar U\otimes U$ and $\bar R_U\colon \un \to U\otimes \bar U$ defined by
\begin{equation}\label{eq:standard-solutions}
R_U(1)=\sum_i\bar\xi_i\otimes\rho_U^{-1/2}\xi_i\ \ \text{and}\ \ \bar R_U(1)=\sum_i\rho_U^{1/2}\xi_i\otimes\bar\xi_i,
\end{equation}
where $\{\xi_i\}_i$ is any orthonormal basis in $H_U$. They solve the conjugate equations for $U$ and $\bar U$, meaning that
$$
(R^*_U\otimes\iota)(\iota\otimes \bar R_U)=\iota_{\bar U}\ \ \text{and}\ \ (\bar R_U^*\otimes\iota)(\iota\otimes R_U)=\iota_U.
$$
The quantum dimension of $U$ is
$$
\dim_qU=\Tr(\rho_U^{\pm1})=\|R_U\|^2=\|\bar R_U\|^2.
$$

\subsection{Quantum group actions}

A (continuous) left action of a compact quantum group $G$ on a C$^*$-algebra $A$ is an injective $*$-homomorphism $\alpha\colon A\to C(G)\otimes A$ such that $(\Delta\otimes\iota)\alpha=(\iota\otimes\alpha)\alpha$ and such that the space $(C(G)\otimes1)\alpha(A)$ is norm dense in $C(G)\otimes A$ (the \emph{Podle\'s condition}). Given such an action, we also say that $A$ is a $G$-C$^*$-algebra.

An element $a\in A$ is called \emph{regular}, if $\alpha(a)$ lies in the algebraic tensor product $\C[G]\otimes_\algg A$. The injectivity of $\alpha$ implies that $(\eps\otimes\iota)\alpha(a)=a$ for regular $a$, where  $\eps$ is the counit of $(\C[G],\Delta)$. It follows that the set of regular elements forms a $*$-subalgebra $\A\subset A$ and $\alpha$ defines a coaction of $(\C[G],\Delta)$ on $\A$.
We have a right $c_c(\hat G)$-module structure on $A$ defined by
$$
\blacktriangleleft\colon A\otimes c_c(\hat G)\to A, \quad a\blacktriangleleft\omega=(\omega\otimes\iota)\alpha(a).
$$
Here we use that every element $\omega\in c_c(\hat G)\subset\C[G]^*$ has the form $\omega=h(\cdot\,x)$ for some $x\in\C[G]$ and therefore extends to a bounded linear functional on $C(G)$, and even to a normal linear functional on $L^\infty(G)=\pi_h(C(G))''$. We have $\A=A\blacktriangleleft c_c(\hat G)$. The Podle\'s condition implies then that $\A$ is dense in~$A$. The converse is also true, see, e.g., \cite{MR2034922}*{Corollary~1.4}.

\smallskip

A left action of $G$ on a von Neumann algebra $N$ is an injective normal unital $*$-ho\-mo\-mor\-phism $\alpha\colon N\to L^\infty(G)\bar\otimes N$ such that $(\Delta\otimes\iota)\alpha=(\iota\otimes\alpha)\alpha$. We then say that $N$ is a $G$-von Neumann algebra. We can define in the same way as above the subalgebra $\NN\subset N$ of regular elements. As a consequence of the Takesaki duality (see~\cite{MR1814995}), the algebra $\NN$ is dense in $N$ in the ultrastrong operator topology. Denote by $\RR(N)$ the norm closure of $\NN$ in $N$. Then the restriction of $\alpha$ to $\RR(N)$ defines an action of $G$ on $\RR(N)$. Thus, we get a functor
$$
\RR\colon (G\text{-von Neumann algebras})\to (G\text{-C}^*\text{-algebras}).
$$

Given two $G$-C$^*$-algebras $A_1$ and $A_2$, we say that a bounded linear map $T\colon A_1\to A_2$ is $G$-equivariant, if it is a $c_c(\hat G)$-module map. We use the same definition for $G$-von Neumann algebras. Furthermore, we can and will consider the situations when one of the algebras is a $G$-C$^*$-algebra and the other is a $G$-von Neumann algebra, the equivariance understood this way still makes sense.

\smallskip

A right action of the dual discrete quantum group $\hat G$ on a C$^*$-algebra $A$ is an injective nondegenerate $*$-homomorphism
$\beta\colon A\to M(A\otimes c_0(\hat G))$
%=\ell^\infty\text{-}\bigoplus_{s\in\Irr(G)}A\otimes B(H_s)
such that $(\iota\otimes\Dhat)\beta=(\beta\otimes\iota)\beta$. We then say that $A$ is a $\hat G$-C$^*$-algebra. The Podle\'s condition -- density of $(1\otimes c_0(\hat G))\beta(A)$ in $A\otimes c_0(\hat G)$ -- is automatically satisfied in this case. Indeed, since $(\iota\otimes\hat\eps)\beta(a)=a$ for all $a\in A$ by the injectivity of $\beta$, this follows, e.g., from \cite{MR2034922}*{Corollary~1.4}.

Given a right action $\beta$ of $\hat G$, by duality we get a left $\C[G]$-module structure on $A$ defined by
\begin{equation}\label{eq:action}
\rhd\colon \C[G]\otimes A\to A,\quad x\rhd a=(\iota\otimes x)\beta(a).
\end{equation}
Then $A$ becomes a $\C[G]$-module algebra and
\begin{equation}\label{estar}
x\rhd a^*=(S(x)^*\rhd a)^*\quad\text{for all}\quad x\in\C[G],\ a\in A.
\end{equation}

Conversely, if a C$^*$-algebra $A$ is a $\C[G]$-module algebra and condition \eqref{estar} is satisfied, then there is a $*$-homomorphism $\beta\colon A\to \prod_{s\in\Irr(G)}A\otimes B(H_s)$ uniquely determined by~\eqref{eq:action}, where $\iota\otimes x$ is well-defined, as every $x\in\C[G]\subset\U(G)^*$ factors through a finite direct sum of the algebras $B(H_s)$. We then have
$$
(\iota\otimes\Dhat)\beta=(\beta\otimes\iota)\beta\colon A\to\prod_{s,t\in\Irr(G)}A\otimes B(H_s)\otimes B(H_t),
$$
$(\iota\otimes\hat\eps)\beta=\iota$ and $(1\otimes c_c(\hat G))\beta(A)=A\otimes_\algg c_c(\hat G)$ (see, e.g., the proof of~\cite{MR2034922}*{Proposition~1.3}). As the image of $\beta$ is automatically contained in
$$
\ell^\infty\text{-}\bigoplus_{s\in\Irr(G)}A\otimes B(H_s)\subset M(A\otimes c_0(\hat G)),
$$
it follows that $\beta$ defines an action of $\hat G$ on $A$. Therefore the actions of discrete quantum groups admit a purely algebraic description. %cf.~\cite{MR3291643}*{Section~1.5}.

A $\hat G$-von Neumann algebra is a $\hat G$-C$^*$-algebra $N$ such that $N$ is a von Neumann algebra and the action map $N\to M(N\otimes c_0(\hat G))=N\bar\otimes\ell^\infty(\hat G)$ is normal. The last condition is equivalent to saying that $\C[G]$ acts on $N$ by normal operators; it doesn't seem to be known whether this condition is really necessary.

\smallskip

Given two $\hat G$-C$^*$-algebras $A_1$ and $A_2$, we say that a linear map $A_1\to A_2$ is $\hat G$-equivariant if it is a $\C[G]$-module map.

\subsection{Yetter--Drinfeld algebras}\label{ssec:YD}

Assume $G$ is a compact quantum group, $A$ is a C$^*$-algebra and we have a left action $\alpha\colon A\to C(G)\otimes A$ of $G$ and a right action $\beta\colon A\to M(A\otimes c_0(\hat G))$ of the dual discrete quantum group~$\hat G$. The Yetter--Drinfeld condition is commutativity of the following diagram:
\begin{center}
	\begin{tikzcd}[scale = 10em]
C(G)\otimes A \arrow[d, "\iota\otimes \beta", swap] & \arrow[l, "\alpha", swap] A \arrow[r, "\beta"] & M(A \otimes c_0(\hat{G}) ) \arrow[d, "\alpha\otimes\iota"] \\
M(C(G)\otimes A\otimes c_0(\hat{G})) \arrow[rr, "\Ad W_{13}", swap] & & M(C(G)\otimes A\otimes c_0(\hat{G}) )
	\end{tikzcd}
\end{center}
If it is satisfied, we say that $A$ is a Yetter--Drinfeld $G$-C$^*$-algebra. (We could equally well say that it is a Yetter--Drinfeld $\hat G$-C$^*$-algebra.) In a similar way we can define Yetter--Drinfeld $G$-von Neumann algebras.

Let $\A\subset A$ be the subalgebra of regular elements (with respect to the $G$-action). Then the Yetter--Drinfeld condition is equivalent to
\begin{equation} \label{eYD1}
\alpha(x\rhd a) =x_{(1)} a_{(1)}S(x_{(3)})\otimes ( x_{(2)}\rhd a_{(2)})
\end{equation}
for all $x\in\C[G]$ and $a\in\A$, where we use Sweedler's sumless notation, so we write $\Delta(x)=x_{(1)}\otimes x_{(2)}$ and $\alpha(a)=a_{(1)}\otimes a_{(2)}$.
This implies that if $A$ is a Yetter--Drinfeld $G$-C$^*$- or $G$-von Neumann algebra, then $\A\subset A$ is a $\C[G]$-submodule. As a consequence, if $N$ is a Yetter--Drinfeld $G$-von Neumann algebra, then $\RR(N)$ is a Yetter--Drinfeld $G$-C$^*$-algebra.

The most important examples of Yetter--Drinfeld $G$-C$^*$-algebras for us are $C(G)$, with the $G$-action defined by $\Delta$ and the $\hat G$-action defined by
$$
C(G)\ni x\mapsto W(x\otimes 1)W^*,
$$
and $c_0(\hat G)$, with the $\hat G$-action given by $\Dhat$ and the $G$-action defined by
$$
c_0(\hat G)\ni\omega\mapsto W^*(1\otimes\omega)W.
$$
In the same way we get the Yetter--Drinfeld $G$-von Neumann algebras $L^\infty(G)$ and $\ell^\infty(\hat G)$.

\smallskip

To have the structure of a Yetter--Drinfeld $G$-algebra on a C$^*$-algebra $A$ is the same as having a left action of the Drinfeld double $D(G)$ of $G$ on $A$~\cite{MR2566309}*{Proposition~3.2}. Let us briefly recall this correspondence; it will be used only in Section~\ref{ssec:Poisson}, as most of the time it is more convenient to use the Yetter--Drinfeld condition directly.

The reduced C$^*$-algebra of continuous functions on $D(G)$ and the coproduct on it are defined~by
\begin{equation}\label{eq:DD}
C(D(G))=C(G)\otimes c_0(\hat G),\qquad \Delta_{D(G)}(a\otimes\omega)=W_{32}\Delta(a)_{13}\Dhat(\omega)_{42}W_{32}^*.
\end{equation}
To recognize that this is exactly the same definition as in~\cite{MR2566309}, note that the dual of $(C(G),\Delta)$ within the general theory of locally compact quantum groups is $(c_0(\hat G),\Dhat^{\mathrm{op}})$ rather than $(c_0(\hat G),\Dhat)$.

Now, for any Yetter--Drinfeld $G$-C$^*$-algebra $A$, the corresponding left action
$$
\gamma\colon A\to M(C(D(G))\otimes A)
$$
of $D(G)$ is defined by $\gamma=\big((\iota\otimes\beta)\alpha(\cdot)\big)_{132}$. For example, this way $\Delta_{D(G)}$ corresponds to
\begin{equation}\label{eq:DD-on-itself}
\alpha=\Delta\otimes\iota,\qquad \beta=W_{13}(\iota\otimes\Dhat)(\cdot)W_{13}^*,
\end{equation}
which are the actions of the quantum subgroups $G$ and $\hat G^{\mathrm{op}}$ of $D(G)$ by left translations.

\smallskip

Given two Yetter--Drinfeld $G$-C$^*$- or $G$-von Neumann algebras $A_1$ and $A_2$, we say that a bounded linear map $A_1\to A_2$ is $D(G)$-equivariant if it is both $G$- and $\hat G$-equivariant.

\section{Furstenberg--Hamana boundaries of Drinfeld doubles}\label{sec:FH}

\subsection{Minimal idempotents in convex contraction semigroups}

As discussed in the introduction, Hamana's construction of an injective envelope has been adapted to several different contexts. We ourselves will need two versions of it. The following proposition is an attempt to capture the essence of Hamana's arguments in one general statement.

\begin{proposition}\label{prop:Hamana}
Assume $X$ is a subspace of a dual Banach space $Y^*$ and $S$ is a convex semigroup of contractive linear maps $X\to X$ such that if we consider $S$ as a set of maps $X\to Y^*$, then $S$ is closed in the topology of pointwise weak$^*$ convergence. Then there is an idempotent $\phi_0\in S$ such that
\begin{equation}\label{eq:minimal-idempotent}
\phi_0\psi\phi_0=\phi_0\quad\text{for all}\quad\psi\in S.
\end{equation}
\end{proposition}

\bp
Define a pre-order on $S$ by
$$
\phi\prec\psi\quad\text{iff}\quad\|\phi(x)\|\le\|\psi(x)\|\quad\text{for all}\quad x\in X.
$$
Since by assumption $S$ is compact in the topology of pointwise weak$^*$ convergence, every decreasing chain in $S$ with respect to this pre-order has a lower bound. By Zorn's lemma we conclude that $S$ has a minimal element $\phi_0$.

Now, assume $\psi\in S$ is such that $\psi(\phi_0(X))\subset\phi_0(X)$. We claim that $\psi|_{\phi_0(X)}$ is the identity map. Consider a cluster point $\Psi\in S$ of
the sequence $\{\frac{1}{n}\sum^n_{k=1}\psi^k\}^\infty_{n=1}$ in $S$. Then $\Psi\psi=\Psi$. As~$\phi_0$ is minimal and $\Psi\phi_0\prec\phi_0$, the map $\Psi|_{\phi_0(X)}$ must be isometric. Hence, for every $x\in \phi_0(X)$, we have
$$
\|\psi(x)-x\|=\|\Psi(\psi(x)-x)\|=\|\Psi\psi(x)-\Psi(x)\|=0,
$$
proving our claim.

Applying the claim to $\psi=\phi_0$, we conclude that $\phi_0$ is an idempotent. Applying it to~$\phi_0\psi$, we get that $\phi_0\psi\phi_0=\phi_0$ for all $\psi\in S$.
\ep

We will call any such $\phi_0$ a \emph{minimal idempotent} in $S$.

\begin{remark}
This agrees with minimality with respect to the standard order on the idempotents defined by
$$
\phi_1\le\phi_2\quad\text{iff}\quad \phi_2\phi_1=\phi_1\phi_2=\phi_1,
$$
cf.~\cite{MR2806689}*{Theorem~2.9}. Indeed, if an idempotent $\phi_0$ satisfies~\eqref{eq:minimal-idempotent}, then it is immediate that it is minimal with respect to $\le$. Conversely, assume $\phi_1$ is a minimal idempotent with respect to $\le$. Let $\phi_0$ be an idempotent satisfying~\eqref{eq:minimal-idempotent}. Replacing $\phi_0$ by $\phi_0\phi_1$, we may assume that $\phi_0\phi_1=\phi_0$. Then the idempotent $\phi_1\phi_0$ satisfies $\phi_1\phi_0\le\phi_1$, hence $\phi_1\phi_0=\phi_1$. Then, for every $\psi\in S$, we have
$$
\phi_1\psi\phi_1=\phi_1\phi_0\psi\phi_1\phi_0=\phi_1\phi_0=\phi_1.
$$

For general compact semigroups $S$, however, minimality with respect to $\le$ does not imply~\eqref{eq:minimal-idempotent}, with any nontrivial compact group giving an example.
\end{remark}

\subsection{Furstenberg--Hamana boundary}

Let $G$ be a compact quantum group. Recall that a unital $G$-C$^*$-algebra $A$ is called $G$-injective if, given unital $G$-C$^*$-algebras $B$ and $C$, a completely isometric $G$-equivariant ucp map $\phi\colon B\to C$ and a $G$-equivariant ucp map $\psi\colon B\to A$, there is a $G$-equivariant ucp map $\tilde\psi\colon C\to A$ making the diagram
$$
\xymatrix{C\ar@{-->}[rd]^{\tilde\psi} & \\ B\ar[u]^{\phi}\ar[r]_{\psi} & A}
$$
commutative. In a similar way one defines $D(G)$-injectivity, that is, injectivity for Yetter--Drinfeld $G$-C$^*$-algebras.

\begin{definition}
We say that a unital Yetter--Drinfeld $G$-C$^*$-algebra $A$ is a $D(G)$-\emph{boundary}, or that the action of $D(G)$ on $A$ is a \emph{boundary action}, if for every unital Yetter--Drinfeld $G$-C$^*$-algebra $B$ and every $D(G)$-equivariant ucp map $\psi\colon A\to B$, the map $\psi$ is automatically completely isometric.
A $D(G)$-boundary $A$ is called a \emph{Furstenberg--Hamana boundary} of $D(G)$, if it is in addition a $D(G)$-injective C$^*$-algebra.
\end{definition}

The term \emph{boundary action} is suggested in~\cite{KKSV}. Adapting Hamana's terminology, a $D(G)$-boundary is also called a $D(G)$-essential extension of the trivial Yetter--Drinfeld $G$-C$^*$-algebra $\C$, and a Furstenberg--Hamana boundary is called a $D(G)$-injective envelope of $\C$.

\begin{theorem}\label{thm:FH}
For any compact quantum group $G$, a Furstenberg--Hamana boundary of $D(G)$ exists and is unique up to isomorphism.
\end{theorem}

We denote the Furstenberg--Hamana boundary of $D(G)$ by $C(\partial_{\FH}D(G))$.

\smallskip

The proof of the theorem follows by now standard lines~\citelist{\cite{MR509025}\cite{MR566081}\cite{MR2806689}\cite{KKSV}}, we mainly have to take care of $G$-equivariance.

We will need the following construction. Given a $\hat G$-C$^*$-algebra $A$, with the action of $\hat G$ given by $\beta\colon A\to M(A\otimes c_0(\hat G))$, and a state $\phi$ on $A$, we have a $\hat G$-equivariant cp map
$$
\PP_\phi=(\phi\otimes\iota)\beta\colon A\to\ell^\infty(\hat G),
$$
called a \emph{Poisson integral}. Such maps are often also called Izumi's Poisson integrals, since their usage in noncommutative probability was pioneered by Izumi in~\cite{MR1916370}.  It has been observed in a number of cases that if we also have an action of $G$ on $A$, then the maps $\PP_\phi$ are $G$-equivariant as well~\cite{MR1916370}*{Lemmas 2.2(3), 3.8(2)}. Let us prove a general result of this sort.

Assume $N$ is a $G$-von Neumann algebra, with the action of $G$ given by $\alpha\colon N\to L^\infty(G)\bar\otimes N$. Consider the von Neumann algebra $N\bar\otimes\ell^\infty(\hat G)$. We have a right action of $\hat G$ on it given by $\iota\otimes\Dhat$. Define a left action of $G$ by
$$
\alpha_W\colon N\bar\otimes\ell^\infty(\hat G)\to L^\infty(G)\bar\otimes N\bar\otimes\ell^\infty(\hat G),\quad\alpha_W(x)=W_{13}^*(\alpha\otimes\iota)(x)W_{13}.
$$
Using~\eqref{eq:pentagon} it is easy to check that this way $N\bar\otimes\ell^\infty(\hat G)$ becomes a Yetter--Drinfeld $G$-von Neumann algebra.

Now, given a $\hat G$-C$^*$-algebra $A$ as above and a completely bounded (cb) linear map $\phi\colon A\to N$, define
$$
\PP_\phi\colon A\to N\bar\otimes\ell^\infty(\hat G)\quad\text{by}\quad \PP_\phi(a)=(\phi\otimes\iota)\beta(a).
$$
To make sense of this definition, note that
$$
\beta(A)\subset\ell^\infty\text{-}\bigoplus_{s\in\Irr(G)}A\otimes B(H_s)\qquad\text{and}\qquad N\bar\otimes\ell^\infty(\hat G)=\ell^\infty\text{-}\bigoplus_{s\in\Irr(G)}N\otimes B(H_s),
$$
so by $\phi\otimes\iota$ we simply mean a collection of maps $A\otimes B(H_s)\to N\otimes B(H_s)$.

\begin{proposition}\label{prop:Izumi-equivariance}
For any Yetter--Drinfeld $G$-C$^*$-algebra $A$ and $G$-von Neumann algebra~$N$, the map $\phi\mapsto\PP_\phi$ defines a one-to-one correspondence between the $G$-equivariant cb maps $\phi\colon A\to N$ and the $D(G)$-equivariant cb maps $\PP\colon A\to N\bar\otimes\ell^\infty(\hat G)$, with the inverse given by $\PP\mapsto (\iota\otimes\hat\eps)\PP$.
\end{proposition}

Obviously, the map $\PP_\phi$ is ucp if and only if $\phi$ is ucp, so we also get a correspondence between the ucp maps.

\bp
It is easy to see that the map $\phi\mapsto\PP_\phi$ defines a one-to-one correspondence between the cb maps $\phi\colon A\to N$ and the $\hat G$-equivariant cb maps $\PP\colon A\to N\bar\otimes\ell^\infty(\hat G)$, with the inverse given by $\PP\mapsto (\iota\otimes\hat\eps)\PP$. Therefore we only need to show that $\PP_\phi$ is $G$-equivariant if and only if $\phi$ is $G$-equivariant. The ``only if'' direction is immediate, as the map $\iota\otimes\hat\eps\colon N\bar\otimes\ell^\infty(\hat G)\to N$ is $G$-equivariant.

Assume now that $\phi$ is $G$-equivariant. Denote by $\alpha_A$ the $G$-action on $A$ and by $\alpha$ the $G$-action on $N$. For an element $y$ of $\prod_s A\otimes B(H_s)$ or $\prod_s N\otimes B(H_s)$,
denote by $y_s$ the component of $y$ in $A\otimes B(H_s)$, resp., $N\otimes B(H_s)$. Then the Yetter--Drinfeld condition for $A$ can be written as
$$
\big((\iota\otimes\beta)\alpha_A(a)\big)_s=(U_s)_{31}^*(\alpha_A\otimes\iota)(\beta(a)_s)(U_s)_{31}.
$$
If $a\in A$ is regular, then the above expressions live in $\C[G]\otimes_\algg A\otimes B(H_s)$. This justifies the following computations for such $a$ and all $s\in\Irr(G)$:
\begin{align*}
\big((\iota\otimes\PP_\phi)\alpha_A(a)\big)_s &= (\iota\otimes\phi\otimes\iota)\big(\big((\iota\otimes\beta)\alpha_A(a)\big)_s\big) \\
&=(\iota\otimes\phi\otimes\iota)\big((U_s)_{31}^*(\alpha_A\otimes\iota)(\beta(a)_s)(U_s)_{31}\big) \\
&=(U_s)_{31}^*\big((\iota\otimes\phi)\alpha_A\otimes\iota\big))(\beta(a)_s)(U_s)_{31} \\
&=(U_s)_{31}^*(\alpha\otimes \iota)(\phi\otimes \iota)(\beta(a)_s)(U_s)_{31} \\
&= \alpha_W(\PP_\phi(a))_s,
\end{align*}
where we used the Yetter-Drinfeld condition in the second equality and the equivariance of $\phi$ in the fourth equality. This implies $G$-equivariance of $\PP_\phi$.
\ep

For the proof of Theorem~\ref{thm:FH} we will need only the following corollary. Consider the Yetter--Drinfeld $G$-von Neumann algebra $\ell^\infty(\hat G)$ as in Section~\ref{ssec:YD}, or in other words, as defined above for $N=\C$. Then $\RR(\ell^\infty(\hat G))$ is a Yetter--Drinfeld $G$-C$^*$-algebra.

\begin{corollary}\label{cor:injectivity}
The C$^*$-algebra $\RR(\ell^\infty(\hat G))$ is $D(G)$-injective.
\end{corollary}

\bp
Take a unital Yetter--Drinfeld $G$-C$^*$-algebra $B$. Since every $G$-equivariant bounded linear map $B\to\ell^\infty(\hat G)$ has image in $\RR(\ell^\infty(\hat G))$, by Proposition~\ref{prop:Izumi-equivariance} we conclude that the map $\phi\mapsto\PP_\phi$ defines a one-to-one correspondence between the $G$-invariant states on $B$ and the $D(G)$-equivariant ucp maps $B\to \RR(\ell^\infty(\hat G))$. This implies the corollary, since every $G$-invariant state on $B$ can be extended to a $G$-invariant state on any given Yetter--Drinfeld $G$-C$^*$-algebra~$C$ containing $B$ as an operator subsystem (see Lemma~\ref{lem:Arveson} for a stronger result).
\ep

\bp[Proof of Theorem~\ref{thm:FH}]
Consider the convex semigroup $S$ of $D(G)$-equivariant ucp maps
$$
\RR(\ell^\infty(\hat G))\to\RR(\ell^\infty(\hat G)).
$$
We can equally well consider all $D(G)$-equivariant ucp maps $\RR(\ell^\infty(\hat G))\to\ell^\infty(\hat G)$, since every such map has image in $\RR(\ell^\infty(\hat G))$ by the $G$-equivariance. Since $\C[G]$ and $c_c(\hat G)$ act on $\ell^\infty(\hat G)$ by normal operators, this implies that $S$ is closed in the topology of pointwise ultraweak operator convergence. Hence, by Proposition~\ref{prop:Hamana}, there is a minimal idempotent~$\phi_0$ in~$S$. Consider the C$^*$-algebra $A=\phi_0(\RR(\ell^\infty(\hat G)))$ with the Choi--Effros product $a\cdot b=\phi_0(ab)$. We claim that $A$ equipped with the actions of $G$ and $\hat G$ obtained by restriction from those on $\RR(\ell^\infty(\hat G))$ is a Furstenberg--Hamana boundary of $D(G)$.

\smallskip
\textit{Step 1. $A$ is a $G$-C$^*$-algebra.}
Denote by $\alpha$ the $G$-action on $\RR(\ell^\infty(\hat G))$ and by $\NN\subset\RR(\ell^\infty(\hat G))$ the subalgebra of regular elements. The $G$-equivariance of $\phi_0$ can be written as $(\iota\otimes\phi_0)\alpha=\alpha\phi_0$. Note also that the C$^*$-algebra $C(G)\otimes A$ can be viewed as an operator subsystem of $C(G)\otimes \RR(\ell^\infty(\hat G))$ equipped with the product $x\cdot y = (\iota\otimes\phi_0)(xy)$. Hence $\alpha_A:=\alpha|_A$ is a well-defined ucp map $A\to C(G)\otimes A$ and, for all $a,b\in A$, we have
$$
\alpha(a\cdot b)=\alpha(\phi_0(ab)) = (\iota\otimes\phi_0)\alpha(ab) = (\iota\otimes\phi_0)\big(\alpha(a)\alpha(b)\big) = \alpha(a)\cdot\alpha(b),
$$
so that $\alpha_A$ is a homomorphism.

For any $a\in \RR(\ell^\infty(\hat G))$ and $x \in C(G)$, we have
$$
(x\otimes 1)\cdot\alpha(\phi_0(a))=(x\otimes 1)\alpha(\phi_0(a))=(\iota\otimes\phi_0)\big((x\otimes 1)\alpha(a)\big).
$$
This implies that the Podle\'s condition for $\alpha_A$ is satisfied. This follows also from the density of the subalgebra $\A=A\cap\NN=\phi_0(\NN)$ of regular elements.

\smallskip
\textit{Step 2. $A$ is a $\hat G$-C$^*$-algebra.} By the $\hat G$-equivariance of $\phi_0$, the space $A$ is a $\C[G]$-submodule of $\RR(\ell^\infty(\hat G))$. Since condition~\eqref{estar} is obviously satisfied for $A$, in order to define an action of~$\hat G$ on~$A$ we only need to check that the module structure respects the product, that is,
$$
x\rhd(a\cdot b)=(x_{(1)}\rhd a)\cdot (x_{(2)}\rhd b)
$$
for all $x\in\C[G]$ and $a,b\in A$. But this follows immediately by applying $\phi_0$ to the same identity for the original product on $\RR(\ell^\infty(\hat{G}))$.

\smallskip
\textit{Step 3. The Yetter--Drinfeld condition is satisfied.} As condition~\eqref{eYD1} holds for the elements of $\NN$, it obviously holds for the elements of $\A\subset\NN$.

\smallskip
\textit{Step 4. $A$ is $D(G)$-injective.} Since $A=\phi_0(\RR(\ell^\infty(\hat G)))$ and $\phi_0$ is an idempotent, this follows from Corollary~\ref{cor:injectivity}.

\smallskip
\textit{Step 5. $A$ is a $D(G)$-boundary.} Assume $\psi\colon A\to B$ is a $D(G)$-equivariant ucp map for some unital Yetter--Drinfeld $G$-C$^*$-algebra~$B$. By the $D(G)$-injectivity of $A$, there is a $D(G)$-equivariant ucp map $\phi\colon B\to A$. Then $\phi\psi\phi_0$ is an element of the semigroup $S$. By the minimality of $\phi_0$ we get that $\phi\psi\phi_0=\phi_0\phi\psi\phi_0=\phi_0$, that is, $\phi\psi$ is the identity map on $A$. Hence~$\psi$ is completely isometric.

\smallskip

It remains to prove the uniqueness up to isomorphism. This is a standard argument, which we recall for the reader's convenience, cf.~\cite{MR566081}*{Theorem~4.1}. Assume $A_1$ and $A_2$ are two Furstenberg--Hamana boundaries of $D(G)$. Then by the injectivity of $A_2$ there is an equivariant ucp map $\gamma\colon A_1\to A_2$. It must be completely isometric, as $A_1$ is a boundary. Hence, by the injectivity of $A_1$, there is an equivariant ucp map $\tilde\gamma\colon A_2\to A_1$ such that $\tilde\gamma\gamma=\iota$. But then $\gamma\tilde\gamma\colon A_2\to A_2$ is an idempotent, which must be completely isometric, as $A_2$ is a boundary. It follows that $\gamma\tilde\gamma=\iota$. Therefore $\gamma$ is a complete order isomorphism of $A_1$ onto $A_2$, hence it is an isomorphism of C$^*$-algebras.
\ep

The following rigidity property of the Furstenberg--Hamana boundary can be deduced from the definition (see~\cite{MR566081}*{Lemma~3.7}), but it is immediate from the construction.

\begin{corollary}\label{cor:rigidity}
The only $D(G)$-equivariant ucp map $C(\partial_{\FH}D(G))\to C(\partial_{\FH}D(G))$ is the identity map.
\end{corollary}

\bp
Using the notation from the proof of Theorem~\ref{thm:FH}, if $\psi\colon A\to A$ is a $D(G)$-equivariant ucp map, then $\psi\phi_0=\phi_0\psi\phi_0=\phi_0$ by the minimality of $\phi_0$, that is, $\psi$ is the identity map.
\ep

The proof of the theorem implies one more nonobvious property of $C(\partial_{\FH}D(G))$. Recall that a Yetter--Drinfeld $G$-C$^*$-algebra $A$ is called \emph{braided-commutative}, if
\begin{equation}\label{eq:br-com}
ab = b_{(2)}(S^{-1}(b_{(1)})\rhd a)
\end{equation}
for all $a,b\in\A$. Both $C(G)$ and $c_0(\hat G)$ are braided-commutative.

\begin{corollary}\label{cor:br-com}
The Yetter--Drinfeld $G$-C$^*$-algebra $C(\partial_{\FH}D(G))$ is braided-commutative.
\end{corollary}

\bp
Using the notation from the proof of Theorem~\ref{thm:FH}, identity~\eqref{eq:br-com} holds for all $a,b\in\A=A\cap\NN$ for the original product on $\NN$. By applying $\phi_0$ we conclude that it also holds for the Choi--Effros product.
\ep

\subsection{Boundary actions}

In this subsection we follow closely the last part of \cite{KKSV}*{Section~4}. Our goal is to prove the following result.

\begin{theorem}\label{thm:boundaries}
Let $G$ be a compact quantum group. Then up to isomorphism the $D(G)$-boundaries are precisely the unital Yetter--Drinfeld $G$-C$^*$-subalgebras of $C(\partial_{\FH}D(G))$. Furthermore, given two such subalgebras $A_1, A_2\subset C(\partial_{\FH}D(G))$, the embedding map $A_1\hookrightarrow A_2$ is the only $D(G)$-equivariant ucp map $A_1\to A_2$ if $A_1\subset A_2$, and there are no such maps if $A_1\not\subset A_2$.
\end{theorem}

Therefore there is a one-to-one correspondence between the isomorphism classes of $D(G)$-boundaries and the unital Yetter--Drinfeld $G$-C$^*$-subalgebras of $C(\partial_{\FH}D(G))$.

\smallskip

For the proof of the theorem we need the following equivariant version of Arveson's extension theorem.

\begin{lemma}\label{lem:Arveson}
Given a unitary representation $U\in M(K(H)\otimes C(G))$ of $G$ on a Hilbert space~$H$, consider the $G$-von Neumann algebra $B(H)$, with the $G$-action
$$
\alpha_U\colon B(H)\to L^\infty(G)\bar\otimes B(H)\quad\text{given by}\quad \alpha_U(T)=U^*_{21}(1\otimes T)U_{21}.
$$
Then $\RR(B(H))$ is a $G$-injective C$^*$-algebra.
\end{lemma}

\bp
Assume $B$ and $C$ are unital $G$-C$^*$-algebras, with the actions of $G$ denoted by $\alpha_B$ and~$\alpha_C$, $\phi\colon B\to C$ is a completely isometric $G$-equivariant ucp map, and $\psi\colon B\to B(H)$ is a $G$-equivariant ucp map. By Arveson's extension theorem, there is a ucp map $\eta\colon C\to B(H)$ such that $\eta\phi=\psi$. Define
$$
\tilde\psi\colon C\to B(H)\quad\text{by}\quad \tilde\psi(c)=(h\otimes\iota)(U_{21}(\iota\otimes\eta)\alpha_C(c)U_{21}^*).
$$
A simple computation shows that $\tilde\psi$ is $G$-equivariant and $\tilde\psi\phi=\psi$.
\ep

If we combine this lemma with Proposition~\ref{prop:Izumi-equivariance} for $N=B(H)$, then we can conclude that the C$^*$-algebra $\RR(B(H)\bar\otimes\ell^\infty(\hat G))$ is $D(G)$-injective. This generalizes Corollary~\ref{cor:injectivity}.

\bp[Proof of Theorem~\ref{thm:boundaries}]
Let $A$ be a $D(G)$-boundary. By the injectivity of $C(\partial_\FH D(G))$ we get a $D(G)$-equivariant ucp map $\phi\colon A \to C(\partial_\FH D(G))$, which then must be completely isometric. We claim that it is a homomorphism.

By considering, e.g., a representation of the crossed product $G\ltimes A$, we can find a unitary representation $U$ of $G$ on a Hilbert space $H$ and a representation $\pi\colon A\to B(H)$ that is $G$-equivariant with respect to the action $T\mapsto U^*_{21}(1\otimes T)U_{21}$ of $G$ on $B(H)$. By Lemma~\ref{lem:Arveson}, there is a $G$-equivariant ucp map $\tilde\pi\colon C(\partial_\FH D(G))\to B(H)$ such that $\tilde\pi\phi=\pi$. Consider also any $D(G)$-equivariant ucp map $\psi\colon\RR(B(H)\bar\otimes\ell^\infty(\hat G))\to C(\partial_\FH D(G))$. By Proposition~\ref{prop:Izumi-equivariance} we then get the following commutative diagram of $D(G)$-equivariant ucp maps:
$$
\xymatrix{C(\partial_\FH D(G))\ar[r]^{\PP_{\tilde\pi}\ \ \ \ } & \RR(B(H)\bar\otimes\ell^\infty(\hat G))\ar[r]^{\ \ \ \ \psi} & C(\partial_\FH D(G))\\
A\ar[u]^\phi\ar[ur]_{\PP_\pi\ \ \ \ } & & }
$$

The map $\PP_\pi$ is a homomorphism, since $\pi$ is. By Corollary~\ref{cor:rigidity}, the composition $\psi\PP_{\tilde\pi}$ is the identity map, hence $E=\PP_{\tilde\pi}\psi$ is an idempotent with image $\PP_{\tilde\pi}(C(\partial_\FH D(G)))$. It follows that~$\PP_{\tilde\pi}$ defines a C$^*$-algebra isomorphism of $C(\partial_\FH D(G))$ onto $\PP_{\tilde\pi}(C(\partial_\FH D(G)))$ equipped with the Choi--Effros product $a\cdot b=E(ab)$. Since the Choi--Effros product coincides with the original product on $\PP_\pi(A)\subset \PP_{\tilde\pi}(C(\partial_\FH D(G)))$, we conclude that $\phi$ is indeed a homomorphism.

\smallskip

Next, let $A$ be a unital Yetter--Drinfeld $G$-C$^*$-subalgebra of $C(\partial_{\FH}D(G))$ and $\phi\colon A\to B$ be a $D(G)$-equivariant ucp map for some unital Yetter--Drinfeld $G$-C$^*$-algebra $B$. Take any $D(G)$-equivariant ucp map $\psi\colon B\to C(\partial_{\FH}D(G))$. The map $\psi\phi\colon A\to C(\partial_{\FH}D(G))$ extends to a $D(G)$-equivariant ucp map $C(\partial_{\FH}D(G))\to C(\partial_{\FH}D(G))$, which then must be the identity map. This shows that $\phi$ is completely isometric, proving that $A$ is a $D(G)$-boundary. The argument shows also that in the case $B=C(\partial_{\FH}D(G))$ and $\psi=\iota$, the map $\phi\colon A\to C(\partial_{\FH}D(G))$ must be the embedding map. This implies the last part of the theorem.
\ep

\section{Comparison with Poisson boundaries}\label{sec:PvsF}

\subsection{Noncommutative Poisson boundaries}\label{ssec:Poisson}

Given a von Neumann algebra $M$, by a Markov operator on $M$ one means a normal ucp map $P\colon M\to M$. The corresponding Poisson boundary~\cite{MR1916370} is defined as the space
$$
H^\infty(M,P)=\{x\in M\mid P(x)=x\}
$$
of $P$-harmonic elements. This space is ultraweakly operator closed in $M$ and it is the image of a ucp projection $M\to H^\infty(M,P)$ obtained as a cluster point of $\{\frac{1}{n}\sum^n_{k=1}P^k\}^\infty_{n=1}$. As a consequence it has a unique structure of a von Neumann algebra such that the embedding $H^\infty(M,P)\to M$ is a normal completely isometric ucp map.

\smallskip

Let $G$ be a compact quantum group. For a normal state $\phi$ on $\ell^\infty(\hat{G})$, consider the convolution operator
$$
P_{\phi} = (\phi \otimes \iota)\Dhat\colon\ell^\infty(\hat G)\to \ell^\infty(\hat G).
$$
We will be mainly interested in the left $G$-invariant normal states $\phi$. Every such state is determined by its restriction to $\ell^\infty(\hat G)^G=Z(\ell^\infty(\hat G))\cong\ell^\infty(\Irr(G))$, that is, by a probability measure $\mu$ on $\Irr(G)$. Explicitly, the $G$-invariant normal state $\phi_\mu$ corresponding to $\mu$ is given by
$\sum_{s\in\Irr(G)}\mu(s)\phi_s$, where $\phi_s$ is the state on $B(H_s)$ defined by
$$
\phi_s=(\dim_q U_s)^{-1}\Tr(\cdot\,\pi_s(\rho)^{-1}).
$$
We will write $P_s$ and $P_\mu$ instead of $P_{\phi_s}$ and $P_{\phi_\mu}$. We will also use the lighter notation $H^\infty(\hat G,\mu)$ for $H^\infty(\ell^\infty(\hat G),P_{\phi_\mu})$.

The Markov operators $P_\mu$ are $G$-equivariant by Proposition~\ref{prop:Izumi-equivariance} (or by \cite{MR1916370}*{Lemma~2.2}). In particular, they leave $\ell^\infty(\Irr(G))$ invariant and therefore define random walks on $\Irr(G)$. In fact, the states~$\phi_\mu$ are the only ones with this property: if $P_\phi$, for a normal state $\phi$, leaves $\ell^\infty(\Irr(G))$ invariant, then $\phi=\phi_\mu$ for some $\mu$, see~\cite{MR2034922}*{Proposition~2.1}. The operators $P_\mu$ are right $\hat G$-equivariant as well, and as a result $H^\infty(\hat G,\mu)$ is a Yetter--Drinfeld $G$-von Neumann algebra.

\smallskip

Random walks can also be considered on $D(G)$, but as the following result shows, for a natural class of states the corresponding Poisson boundary is described entirely in terms of $\hat G$. %This is not going to be used later, but explains why we concentrate on random walks on~$\hat G$.

For a probability measure $\mu$ on $\Irr(G)$, consider the normal state $h\otimes\phi_\mu$ on $L^\infty(D(G))=L^\infty(G)\bar\otimes\ell^\infty(\hat G)$ and the Markov operator $(\iota\otimes(h\otimes\phi_\mu))\Delta_{D(G)}$ on $L^\infty(D(G))$. We denote by $H^\infty(D(G),\mu)$ the corresponding Poisson boundary.

\begin{proposition}\label{poisson-boundary-D(G)}
For any compact quantum group $G$ and any probability measure $\mu$ on $\Irr(G)$, we have an isomorphism
$$
H^\infty(\hat G,\mu)\cong H^\infty(D(G),\mu),\quad \omega\mapsto W^*(1\otimes \omega)W,
$$
of Yetter--Drinfeld $G$-von Neumann algebras.
\end{proposition}

\bp
Denote by $Q$ the Markov operator $(\iota\otimes(h\otimes\phi_\mu))\Delta_{D(G)}$ on $L^\infty(D(G))$. By the definition~\eqref{eq:DD} of $\Delta_{D(G)}$ we have
$$
Q(a\otimes\omega)=(\iota\otimes h\otimes\iota)\big(W_{23}(\Delta(a)\otimes P_\mu(\omega))W_{23}^*\big).
$$
Using first that $W_{23}=W_{13}^*(\Delta\otimes\iota)(W)$ by~\eqref{eq:pentagon} and then the invariance of the Haar state $h$, we see that this equals
$$
W^*(\iota\otimes h\otimes\iota)(\Delta\otimes\iota)\big(W(a\otimes  P_\mu(\omega))W^*\big)W
=W^*(h(\cdot)1\otimes\iota)\big(W(a\otimes  P_\mu(\omega))W^*\big)W.
$$
It follows that
$$
H^\infty(D(G),\mu)\subset\operatorname{Im}Q\subset W^*(1\otimes\ell^\infty(\hat G))W.
$$

Next, consider an element of the form $x=W^*(1\otimes \omega)W$, $\omega\in\ell^\infty(\hat G)$. Then the above computation gives
$$
Q(x)=W^*(h(\cdot)1\otimes\iota)\big(W(\iota\otimes  P_\mu)\big(W^*(1\otimes\omega)W\big)W^*\big)W.
$$
As $(\iota\otimes  P_\mu)\big(W^*(1\otimes\omega)W\big)=W^*(1\otimes P_\mu(\omega))W$ by the $G$-equivariance of $P_\mu$, we thus get
$$
Q(x)=W^*(1\otimes P_\mu(\omega))W,
$$
so that $x$ is $Q$-harmonic if and only if $\omega$ in $P_\mu$-harmonic. This proves the isomorphism $H^\infty(\hat G,\mu)\cong H^\infty(D(G),\mu)$ stated in the proposition. The claim that this isomorphism respects the actions of $G$ and $\hat G$ is easy to check using~\eqref{eq:pentagon}, if we recall that the actions on $H^\infty(D(G),\mu)\subset L^\infty(D(G))$ are given by~\eqref{eq:DD-on-itself}.
\ep

\begin{remark}
In the Kac case, and only in this case, the states of the form $h\otimes\phi_\mu$ are exactly the normal states on $L^\infty(D(G))$ that are invariant under the actions of the quantum subgroup~$G$ of $D(G)$ by left and right translations.
\end{remark}

\subsection{Quantum groups with weakly amenable dimension functions}\label{ssec:weak-amenable}

Recall that given a right action $\beta\colon A\to M(A\otimes c_0(\hat G))$ of $\hat G$ and a normal state $\phi$ on $\ell^\infty(\hat G)$, a state~$\omega$ on $A$ is called \emph{$\phi$-stationary} if $\omega*\phi:=(\omega\otimes\phi)\beta=\omega$, or in other words, if $\phi\PP_\omega=\omega$~\citelist{\cite{MR146298}\cite{KKSV}}. It is easy to see that the Poisson integral $\PP_\omega\colon A\to \ell^\infty(\hat G)$ has image in $H^\infty(\ell^\infty(\hat G),P_\phi)$ if and only if $\omega$ is $\phi$-stationary.

\begin{proposition}\label{prop:PvsFH}
Assume $G$ is a compact quantum group and $\mu$ is a probability measure on~$\Irr(G)$. Then there is a $D(G)$-equivariant ucp map $\psi\colon\RR(H^\infty(\hat G,\mu))\to C(\partial_\FH D(G))$. Any such map $\psi$ is surjective and the following conditions are equivalent:
\begin{enumerate}
  \item[(1)] $\psi$ is an isomorphism;
  \item[(2)] the only $D(G)$-equivariant ucp map $\RR(H^\infty(\hat G,\mu))\to\RR(H^\infty(\hat G,\mu))$ is the identity map;
  \item[(3)] the state $\hat\eps|_{\RR(H^\infty(\hat G,\mu))}$ is the only $G$-invariant $\phi_\mu$-stationary state on $\RR(H^\infty(\hat G,\mu))$;
  \item[(4)] $\RR(H^\infty(\hat G,\mu))$ is a $D(G)$-boundary.
\end{enumerate}
\end{proposition}

\bp
A $D(G)$-equivariant ucp map $\psi\colon\RR(H^\infty(\hat G,\mu))\to C(\partial_\FH D(G))$ exists by the $D(G)$-injectivity of $C(\partial_\FH D(G))$. On the other hand, $\RR(H^\infty(\hat G,\mu))$ is $D(G)$-injective as well, being the image of an equivariant ucp projection $\RR(\ell^\infty(\hat G))\to \RR(H^\infty(\hat G,\mu))$. Hence there is an equivariant ucp map $\eta\colon C(\partial_\FH D(G))\to \RR(H^\infty(\hat G,\mu))$. Then $\psi\eta=\iota$ by Corollary~\ref{cor:rigidity}, which shows that $\psi$ is surjective. Let us show equivalence of (1)--(4).

\smallskip

The implication (1)$\Rightarrow$(2) follows again from Corollary~\ref{cor:rigidity}.

\smallskip

Assuming (2), let $\omega$ be a $G$-invariant $\phi_\mu$-stationary state. Then $\PP_\omega\colon \RR(H^\infty(\hat G,\mu))\to\ell^\infty(\hat G)$ has image in $H^\infty(\hat G,\mu)$, hence in $\RR(H^\infty(\hat G,\mu))$ by the $G$-equivariance. But then by assumption it must be the identity map. Hence $\hat\eps=\hat\eps\PP_\omega=\omega$ on $\RR(H^\infty(\hat G,\mu))$. Note that the state $\hat\eps|_{\RR(H^\infty(\hat G,\mu))}$ is indeed $\phi_\mu$-stationary, since $\hat\eps=\phi_\mu$ on $H^\infty(\hat G,\mu))$. Thus, (2)$\Rightarrow$(3).

\smallskip

Next, assume that (3) holds and we are given a $D(G)$-equivariant ucp map $\psi\colon\RR(H^\infty(\hat G,\mu))\to A$ for a unital Yetter--Drinfeld $G$-C$^*$-algebra $A$. As $\RR(H^\infty(\hat G,\mu))$ is $D(G)$-injective, there is also an equivariant ucp $\eta\colon A\to \RR(H^\infty(\hat G,\mu))$. Consider the state $\omega=\hat\eps\eta\psi$. Then, by Proposition~\ref{prop:Izumi-equivariance}, we have $\eta\psi=\PP_\omega$. But since $\hat\eps|_{\RR(H^\infty(\hat G,\mu))}$ is $\phi_\mu$-stationary, the state~$\omega$ is $\phi_\mu$-stationary as well, hence $\omega=\hat\eps|_{\RR(H^\infty(\hat G,\mu))}$ and therefore $\PP_\omega$ is the identity map. Thus, $\eta\psi=\iota$, which implies that $\psi$ is a complete isometry. This shows that (3)$\Rightarrow$(4).

\smallskip

The implication (4)$\Rightarrow$(1) follows, e.g., from the surjectivity of $\psi$ and Theorem~\ref{thm:boundaries}.
\ep

\begin{remark}\label{rem:stationary}
The proof of the implication (3)$\Rightarrow$(4) shows that, for any $G$ and $\mu$, if $B\subset \RR(H^\infty(\hat G,\mu))$ is a unital Yetter--Drinfeld $G$-C$^*$-subalgebra such that $\hat\eps|_B$ is the only $G$-invariant $\phi_\mu$-stationary state on $B$, then $B$ is a $D(G)$-boundary, cf.~\cite{KKSV}*{Theorem~4.19}.\hfill$\diamondsuit$
\end{remark}

The simplest case when condition (3) in Proposition~\ref{prop:PvsFH} is satisfied is when $G$ acts ergodically on $H^\infty(\hat G,\mu)$ and therefore $\RR(H^\infty(\hat G,\mu))$ has a unique $G$-invariant state. This is probably the closest to the original setup and motivation of Furstenberg, see Corollary to Theorem~3.1 in~\cite{MR146298}. It is very much possible that Proposition~\ref{prop:PvsFH} cannot be applied in any other case: if $H^\infty(\hat G,\mu)^G$ is nontrivial, the C$^*$-algebra $\RR(H^\infty(\hat G,\mu))$ might be too large to be the Furstenberg--Hamana boundary.

The condition $H^\infty(\hat G,\mu)^G=\C1$ is equivalent to triviality of the Poisson boundary of the classical random walk on $\Irr(G)$ defined by $P_\mu|_{\ell^\infty(\Irr(G))}$. A probability measure $\mu$ satisfying this condition is called \emph{ergodic}. As a byproduct of the isomorphism $C(\partial_\FH D(G)))\cong\RR(H^\infty(\hat G,\mu))$, we conclude that up to isomorphism the Yetter--Drinfeld $G$-von Neumann algebra $H^\infty(\hat G,\mu)$ is independent of an ergodic measure $\mu$. In fact, a stronger result is known and is easy to prove. Define
$$
H^\infty(\hat G)=\{x\in \ell^\infty(\hat G)\mid P_s(x)=x\ \text{for all}\ s\in\Irr(G)\}.
$$
The elements of $H^\infty(\hat G)$ are called \emph{absolutely harmonic}~\cite{MR3782061}.

\begin{lemma}
For any ergodic probability measure $\mu$ on $\Irr(G)$, we have $H^\infty(\hat G,\mu)=H^\infty(\hat G)$.
\end{lemma}

\bp
Take a probability measure $\nu$ on $\Irr(G)$. Since the action of $G$ on $H^\infty(\hat G,\mu)$ is ergodic, the state $\hat\eps$ is the only normal $G$-invariant state on $H^\infty(\hat G,\mu)$. Hence $\phi_\nu=\hat\eps$ on $H^\infty(\hat G,\mu)$, and therefore $P_\nu|_{H^\infty(\hat G,\mu)}$ is the identity map. It follows that  $H^\infty(\hat G,\mu)\subset H^\infty(\hat G)$. The opposite inclusion is obvious.
\ep

Recall that the quantum dimension function $U\mapsto\dim_qU$ on $\Rep G$ is said to be \emph{weakly amenable}, if there is a state on $\ell^\infty(\Irr(G))=Z(\ell^\infty(\hat G))$ that is invariant under $P_s$ for all $s\in\Irr(G)$. By a version of the Furstenberg--Kaimanovich--Vershik--Rosenblatt theorem~\cite{MR1749868}*{Theorem~2.5}, the quantum groups admitting ergodic measures are exactly the ones that have weakly amenable quantum dimension functions and countable $\Irr(G)$. The corresponding Poisson boundaries have been shown to have a universal property~\citelist{\cite{MR3291643}\cite{MR3679617}}. We want to show next that this property implies that $H^\infty(\hat G)$ leads to a model of the Furstenberg--Hamana boundary of~$D(G)$ without the countability assumption on $\Irr(G)$. We need some preparation to formulate the result.

\smallskip

Recall (see \cite{MR3679617}*{Section~2} and references there) that a dimension function $d$ on $\Rep G$ is called \emph{amenable}, if
$$
d(U)=\|\Gamma_U\|
$$
for all representations $U\in\Rep G$, where $\Gamma_U\in B(\ell^2(\Irr(G)))$ is the matrix defined by
$$
\Gamma_U=(\dim\Mor(U_s,U\otimes U_t))_{s,t}.
$$
Since the inequality $\|\Gamma_U\|\le d(U)$ holds for any dimension function, an amenable dimension function is unique if it exists, and when it exists, it is the smallest dimension function on~$\Rep G$. Coamenability of $G$ is equivalent to amenability of the classical dimension function $U\mapsto\dim H_U$ on $\Rep G$. In particular, the quantum dimension function is amenable if and only if $G$ is coamenable and of Kac type.

Now, consider a unital braided-commutative Yetter--Drinfeld C$^*$-algebra $A$, with coactions
$$
\alpha\colon A\to C(G)\otimes A\qquad\text{and}\qquad \beta\colon A\to M(A\otimes c_0(\hat G)).
$$
Then by~\cite{MR3291643} we can associate to $A$ a C$^*$-tensor category $\CC_A$ together a unitary tensor functor $\Rep G\to\CC_A$. A convenient way of doing this is to start with $\Rep G$, enlarge the morphism spaces and then consider the idempotent completion of the new category. Namely, the new morphism spaces are defined by
$$
\CC_A(U,V)=\{T\in A\otimes B(H_U,H_V): V_{31}^*(\alpha\otimes\iota)(T)U_{31}=1\otimes T\}.
$$
We will write $\CC_A(U)$ for $\CC_A(U,U)=\End_{\CC_A}(U)$. The tensor products of morphisms are described by the following rules:
\begin{equation}\label{eq:tensor-structure}
T\otimes\iota_Y=T\otimes1\in A\otimes B(H_U,H_V)\otimes B(H_Y),\qquad \iota_Y\otimes T=(\beta_Y\otimes\iota)(T)
\end{equation}
for all $Y\in\Rep G$, where $\beta_Y=(\iota\otimes\pi_Y)\beta\colon A\to A\otimes B(H_Y)$. The morphisms $T$ in the category $\Rep G$ are viewed as morphisms in $\CC_A$ via the map $T\mapsto 1\otimes T$.

If the action of $G$ on $A$ is ergodic, then $\CC_A$ becomes a rigid C$^*$-tensor category with simple unit, so it has a well-defined intrinsic dimension function. Standard solutions of the conjugate equations in $\CC_A$ can be expressed in terms of the solutions~\eqref{eq:standard-solutions}, which might not be standard in~$\CC_A$. Namely, cf. \cite{MR3679617}*{Lemma~4.7}, for every $U\in\Rep G$ there is a unique positive invertible element $a_U\in\CC_A(U)=(A\otimes B(H_U))^G$ such that
$$
(\iota\otimes a_U^{1/2})R_U,\qquad (a_U^{-1/2}\otimes\iota)\bar R_U
$$
form a standard solution of the conjugation equation for $U$ and $\bar U$ in $\CC_A$. In other words, $a_U>0$ is characterized by the properties that
the scalars
$$
R_U^*(\iota\otimes a_U)R_U\quad\text{and}\quad \bar R_U^*(a_U^{-1}\otimes\iota)\bar R_U
$$
in $\CC_A(\un)=A^G=\C1$ are equal and their product is minimal possible. The dimension $d^{\CC_A}(U)$ of $U\in\Rep G$ in $\CC_A$ is therefore given by
$$
d^{\CC_A}(U)1_A=(\dim_q U)(\iota\otimes\psi_U)(a_U^{-1}),
$$
where
$$
\psi_U=(\dim_q U)^{-1}\Tr(\cdot\,\rho_U).
$$

\begin{example}\label{ex:aU}
Assume $K$ is a closed quantum subgroup of $G$, so that we have a surjective homomorphism $\pi\colon\C[G]\to\C[K]$ of Hopf $*$-algebras. Consider the Yetter--Drinfeld $G$-C$^*$-subalgebra $A=C(G/K)$ of $C(G)$. In this case $\CC_A$ is monoidally equivalent to $\Rep K$, see~\cite{MR3291643}*{Section~3.1}, so the corresponding dimension function on $\Rep G$ is defined by the quantum dimension function for~$K$. If $\rho^G$ and $\rho^K$ denote the Woronowicz characters for $G$ and $K$, then the elements $a_U$ are given~by
$$
a_U=U_{21}(1\otimes\rho^G_U(\rho^K_U)^{-1})U_{21}^*\quad\text{for all}\quad U\in\Rep G.
$$
Note (see~\cite{MR2200270}*{Lemma~2.7}) that $\rho^G_U=\rho^K_UT=T\rho^K_U$ for a positive operator $T\in\End_K(H_U)$, since every right $K$-invariant state on $B(H_U)$ has the form $\Tr(\cdot T\rho^K_U)$ for a positive $T\in\End_K(H_U)$. In particular, the operator $\rho^G_U(\rho^K_U)^{-1}=(\rho^K_U)^{-1}\rho^G_U$ is indeed positive and $\pi$ intertwines the scaling groups of~$G$ and~$K$.\hfill$\diamondsuit$
\end{example}

We are now ready to reformulate, and slightly extend, the main results of~\cite{MR3679617} in the quantum group setting.

\begin{theorem}\label{thm:NY-extension}
Assume $G$ is a compact quantum group with weakly amenable quantum dimension function. Then there is a unique up to isomorphism unital braided-commutative Yetter--Drinfeld G-C$^*$-algebra $A$ with $A^G=\C1$ such that the following conditions are satisfied:
\begin{enumerate}
\item the dimension function on $\Rep G$ defined by $\CC_A$ is amenable, that is,
$$
\|\Gamma_U\|1_A=(\dim_q U)(\iota\otimes\psi_U)(a_U^{-1})\quad\text{for all}\quad U\in\Rep G;
$$
\item if $B$ is another unital braided-commutative Yetter--Drinfeld G-C$^*$-algebra such that $B^G=\C1$ and $B$ defines the amenable dimension function on $\Rep G$, then there is a unique unital $*$-homomorphism $A\to B$ of Yetter--Drinfeld G-C$^*$-algebras.
\end{enumerate}

Furthermore, if $\phi$ is the unique $G$-invariant state on $A$, then the Poisson integral $\PP_\phi\colon A\to\ell^\infty(\hat G)$ is completely isometric and extends to a normal completely isometric map of $\pi_\phi(A)''$ onto $H^\infty(\hat G)$.
\end{theorem}

We denote the Yetter--Drinfeld $G$-C$^*$-algebra $A$ given by Theorem~\ref{thm:NY-extension} by $C(\partial_\Pi\hat G)$ and call~$\partial_\Pi\hat G$ the \emph{Poisson boundary} of~$\hat G$.

\bp
When $\Irr(G)$ is countable, this is essentially Theorems 4.1 and 5.1 in~\cite{MR3679617} combined with the correspondence between braided-commutative Yetter--Drinfeld C$^*$-algebras and tensor functors~\cite{MR3291643}. The only part that requires an additional explanation is (ii), namely, that a homomorphism $A\to B$ is truly unique, not just unique up to an automorphism of $B$. But this follows from the proof of~\cite{MR3679617}*{Theorem~4.1}. Namely, denote by $b_U\in (B\otimes B(H_U))^G$ the elements defined similarly to $a_U$. Then any $D(G)$-equivariant unital $*$-homomorphism $\eta\colon A\to B$ must satisfy
\begin{equation}\label{eq:unique-Poisson}
(\eta\otimes\iota)(a_U)=b_U\quad\text{for all}\quad U\in\Rep G.
\end{equation}
On the other hand, by the proof of~\cite{MR3679617}*{Theorem~4.1}, the morphism spaces $\CC_A(U,V)$ are generated by the morphisms in $\Rep G$, the morphisms $a_U$ and their tensor products. It follows that $\eta$ is completely determined by~\eqref{eq:unique-Poisson}.

\smallskip

To deal with the general case, we will first construct a net of quotients $G_i$ of $G$ with countable $\Irr(G_i)$ and weakly amenable quantum dimension functions.

Since there is an invariant state on $\ell^\infty(\Irr(G))$, a standard argument (see, e.g., \cite{MR0251549}*{\S2.4}) shows that there is a net $(m_j)_j$ of normal states such that $m_jP_s|_{\ell^\infty(\Irr(G))}-m_j\to0$ in norm for all $s\in\Irr(G)$. The states $m_j$, viewed as elements of $\ell^1(\Irr(G))$, have at most countable supports $I_j\subset\Irr(G)$.

Take any countable set $X_0\subset\Irr(G)$. Then we can find a sequence $\{j_n\}^\infty_{n=1}$ such that $m_{j_n}P_s|_{\ell^\infty(\Irr(G))}-m_{j_n}\to0$ for all $s\in X_0$. Let $X_1\subset \Irr(G)$ consist of all $s\in\Irr(G)$ such that $U_s$ is a subrepresentation of a tensor product of representations $U_t$ and their conjugates for $t\in X_0\cup(\cup^\infty_{n=1} I_{j_n})$. The set $X_1$ is countable. Repeat the same procedure with $X_0$ replaced by $X_1$, and so on. We thus get a sequence of countable sets $X_0\subset X_1\subset\dots$. Let $G_0$ be the quotient of~$G$ with $\Irr(G_0)=\cup_{n\ge0}X_n$. In other words, $\C[G_0]\subset\C[G]$ is spanned by the matrix coefficients of~$U_s$ for all $s\in \cup_{n\ge0}X_n$. By construction, we can find a sequence $\{j_n'\}^\infty_{n=1}$ such that $I_{j_n'}\subset\Irr(G_0)$ and $m_{j'_n}P_s|_{\ell^\infty(\Irr(G))}-m_{j'_n}\to0$ for all $s\in\Irr(G_0)$. This implies that if we identify $\ell^\infty(\Irr(G_0))$ with a direct summand of $\ell^\infty(\Irr(G))$ and view the states $m_{j'_n}$ as states on $\ell^\infty(\Irr(G_0))$, then a weak$^*$ cluster point of the sequence $\{m_{j'_n}\}^\infty_{n=1}$ is an invariant mean on $\ell^\infty(\Irr(G_0))$. Therefore $G_0$ has weakly amenable quantum dimension function.

We apply the above procedure to every countable subset $X_0$ of $\Irr(G)$ and this way get a collection of quotients $G_i$ of $G$ with countable $\Irr(G_i)$ and weakly amenable quantum dimension functions. Define a partial order on this collection by the inclusions $\Irr(G_i)\subset\Irr(G_k)$. This gives us the required net of quotients of $G$.

\smallskip

For every $k$, consider the Yetter--Drinfeld $G_k$-C$^*$-algebra $A_k$ given by the theorem in the countable case. Note that if $U\in\Rep G_k$ is viewed as a representation of $G$, then the matrix $(\dim\Mor(U_s,U\otimes U_t))_{s,t\in\Irr(G_k)}$ is only a corner of $\Gamma_U$, but its norm equals that of $\Gamma_U$, see the proof of~\cite{MR1644299}*{Proposition~4.8}. Therefore
\begin{equation}\label{eq:norm-hereditary}
d^{\CC_{A_k}}(U)=\|\Gamma_U\|\quad\text{for all}\quad U\in\Rep G_k.
\end{equation}

Assume now that $\Irr(G_i)\subset\Irr(G_k)$ for some $i$ and $k$. Consider the C$^*$-algebra $A_{ik}\subset A_k$ that is the closure of the span of the spectral subspaces of $A_k$ corresponding to $s\in\Irr(G_i)$. Then $A_{ik}$ has the structure of a braided-commutative Yetter--Drinfeld $G_i$-C$^*$-algebra. For every $U\in\Rep G_i$, we have
$$
(A_k\otimes B(H_U))^{G_k}=(A_{ik}\otimes B(H_U))^{G_i},
$$
where on the left hand side we view $U$ as a representation of $G_k$. This and~\eqref{eq:norm-hereditary} imply that $A_{ik}$ defines the amenable dimension function on $\Rep G_i$. It follows that there exists a unique unital $*$-homomorphism $A_i\to A_{ik}$ of Yetter--Drinfeld $G_i$-C$^*$-algebras. We thus get an inductive system of C$^*$-algebras $A_i$. Let $A$ be the limit of this inductive system.

It follows almost immediately by the construction that $A$ has the structure of a braided-commutative Yetter--Drinfeld $G$-C$^*$-algebra and satisfies properties (i) and (ii). It remains to prove the claim about the Poisson integral.

\smallskip

Consider the von Neumann algebra $M=\pi_\phi(A)''$. The state $\phi$ is faithful on $A$, as $\phi(\cdot)1=(h\otimes\iota)\alpha$. Hence we can view $A$ as a subalgebra of $M$. We continue to denote by $\phi$ the normal state $(\cdot\,\xi_\phi,\xi_\phi)$ on $M$. The quantum group $G$ acts ergodically on $M$ and $\phi$ is the unique normal $G$-invariant state on $M$, hence $\phi$ is faithful on $M$ for the same reason as that $\phi|_A$ is faithful.

For every index $i$, let $M_i=\pi_\phi(A_i)''\subset M$. By the faithfulness of $\phi$ on $M$, we have $M_i\cong\pi_{\phi_i}(A_i)''$, where $\phi_i=\phi|_{A_i}$.
Since the theorem is true for quantum groups with countable isomorphism classes of irreducible representations, we already know that $\PP_{\phi_i}$ extends to an isomorphism of von Neumann algebras $M_i\cong H^\infty(\hat G_i)=H^\infty(\hat G_i,\mu_i)$, where $\mu_i$ is any ergodic measure on $\Irr(G_i)$.

Denote by $\pi_{i}$ the map $\ell^\infty(\hat G)\to\ell^\infty(\hat G_i)$ dual to the embedding $\C[G_i]\to\C[G]$. If $\Irr(G_i)\subset\Irr(G_k)$, define in a similar way $\pi_{ki}\colon \ell^\infty(\hat G_k)\to\ell^\infty(\hat G_i)$. We then have the following commutative diagram:
\begin{equation*}
\xymatrix{M_k\ar[r]^{\PP_{\phi_k}\ \ \ \ } & H^\infty(\hat G_k)\ar[d]^{\pi_{ki}} \\ M_i\ar[u]\ar[r]_{\PP_{\phi_i}\ \ \ \ } & H^\infty(\hat G_i),}
\end{equation*}
where $M_i\to M_k$ is the embedding map. Since the horizontal maps are isomorphisms, it follows that $\pi_{ki}(H^\infty(\hat G_k))=H^\infty(\hat G_i)$ and we have a conditional expectation $E_{ki}\colon M_k\to M_i$ such that $\PP_{\phi_i}E_{ki}=\pi_{ki}\PP_{\phi_k}$. This conditional expectation is $\phi_k$-preserving. The existence of these conditional expectations for all $k$ large enough implies that there is a (necessarily unique) $\phi$-preserving conditional expectation $E_i\colon M\to M_i$. We remark that the existence of $E_i$ follows also from the description of the modular group of $\phi$ given in~\cite{MR2202309}*{Proposition~2.10}.

Now, take $x\in M$. Let $y_i=\PP_{\phi_i}(E_i(x))\in H^\infty(\hat G_i)$. As $\pi_{ki}(y_k)=y_i$ for $\Irr(G_i)\subset\Irr(G_k)$, there is a unique element $y\in\ell^\infty(\hat G)$ such that $y_i=\pi_i(y)$ for all $i$. As the homomorphisms~$\pi_i$ respect the comultiplications, it is easy to see that we must have $y\in H^\infty(\hat G)$. We define $\PP(x)=y$. It is then straightforward to check that $\PP\colon M\to H^\infty(\hat G)$ is a normal ucp map that agrees with $\PP_\phi$ on~$A$. As $E_i(x)\to x$ in the ultrastrong operator topology, we have $\|x_i\|\nearrow\|x\|$. It is also clear that $\|y_i\|\nearrow\|y\|$. It follows that $\|x\|=\|y\|$, so that $\PP$ is an isometric map. The same argument applies to the matrix algebras over $M$ and $\ell^\infty(\hat G)$, hence $\PP$ is completely isometric.

What is left to prove is that $\PP(M)=H^\infty(\hat G)$. Take $y\in H^\infty(\hat G)$. For every~$i$ there is a unique $x_i\in M_i$ such that $\pi_i(y)=\PP_{\phi_i}(x_i)$. Then $\|x_i\|=\|\pi_i(y)\|\le\|y\|$. We also have $E_{ki}(x_k)=x_i$ if $\Irr(G_i)\subset\Irr(G_k)$. It follows that the net $(x_i)_i$ converges in the ultrastrong operator topology to a unique element $x\in M$ such that $E_i(x)=x_i$ for all $i$. Then $\PP(x)=y$.
\ep

The theorem implies that $H^\infty(\hat G)$ has a von Neumann algebra structure. Not surprisingly, it is given by a Choi--Effros product by the following elaboration on the proof.

\begin{corollary}\label{cor:h-infty-projection}
There is a $D(G)$-equivariant ucp projection $\ell^\infty(\hat G)\to H^\infty(\hat G)$.
\end{corollary}

\bp
We will use the net of quotients $G_i$ of $G$ constructed in the proof of Theorem~\ref{thm:NY-extension}. For every $i$, choose a $D(G_i)$-equivariant ucp projection $e_i\colon\ell^\infty(\hat G_i)\to H^\infty(\hat G_i)$. Consider the maps $e_i\pi_i\colon\ell^\infty(\hat G)\to H^\infty(\hat G_i)$. Identifying $\ell^\infty(\hat G_i)$ with a direct summand of~$\ell^\infty(\hat G)$, we can view $e_i\pi_i$ as a cp map on $\ell^\infty(\hat G)$. This map is $G$- and $\hat G_i$-equivariant. It follows that any cluster point $e\colon\ell^\infty(\hat G)\to\ell^\infty(\hat G)$ of the net $(e_i\pi_i)_i$ is a $D(G)$-equivariant cp map.

We need to show that $e$ is a projection onto $H^\infty(\hat G)$. Take $x\in\ell^\infty(\hat G)$. If $\Irr(G_i)\subset\Irr(G_k)$, then $\pi_i(e_k\pi_k(x))\in H^\infty(\hat G_i)$, hence by passing to the limit we get that $\pi_i(e(x))\in H^\infty(\hat G_i)$ for all $i$. As we already used in the proof of Theorem~\ref{thm:NY-extension}, this implies that $e(x)\in H^\infty(\hat G)$. Furthermore, if we start with $x\in H^\infty(\hat G)$, then  $\pi_i(e_k\pi_k(x))=\pi_i(x)$, hence $\pi_i(e(x))=\pi_i(x)$ for all $i$ and therefore $e(x)=x$.
\ep

\begin{corollary}
If $\phi$ is the unique $G$-invariant state on $C(\partial_\Pi\hat G)$, then $\PP_\phi$ is a $D(G)$-equivariant isomorphism of $C(\partial_\Pi\hat G)$ onto $\RR(H^\infty(\hat G))$.
\end{corollary}

\bp
Since $G$ acts ergodically on $H^\infty(\hat G)$, all spectral subspaces of $H^\infty(\hat G)$ are finite dimensional and hence $\RR(H^\infty(\hat G))$ is the unique ultrastrongly operator dense $G$-C$^*$-subalgebra of~$H^\infty(\hat G)$. It follows that $\PP_\phi(C(\partial_\Pi\hat G))=\RR(H^\infty(\hat G))$.
\ep

We now return to the Furstenberg--Hamana boundaries.

\begin{theorem} \label{thm:FH-weak-amenable}
For any compact quantum group $G$ with weakly amenable quantum dimension function, we have an isomorphism
$$
C(\partial_\FH D(G))\cong C(\partial_\Pi\hat G)
$$
of Yetter--Drinfeld $G$-C$^*$-algebras.
\end{theorem}

\bp
When $\Irr(G)$ is countable, then $C(\partial_\Pi \hat G)\cong\RR(H^\infty(\hat G,\mu))$ for any ergodic probability measure $\mu$ and the theorem follows from Proposition~\ref{prop:PvsFH}. In the general case the arguments are similar.

The C$^*$-algebra $C(\partial_\Pi\hat G)\cong \RR(H^\infty(\hat G))$ is $D(G)$-injective by Corollary~\ref{cor:h-infty-projection}. To show that it is a $D(G)$-boundary, assume $B$ is a unital Yetter--Drinfeld $G$-C$^*$-algebra and $\psi\colon C(\partial_\Pi\hat G)\to B$ is a $D(G)$-equivariant ucp map. Take a $G$-invariant state~$\omega$ on~$B$. Then $\omega\psi=\phi$, the unique $G$-invariant state on $C(\partial_\Pi\hat G)$. Hence $\PP_{\omega}\psi=\PP_\phi$ is completely isometric by Theorem~\ref{thm:NY-extension}. It follows that $\psi$ is completely isometric as well.
\ep

\begin{corollary}\label{cor:FH-G-ergodic}
For any compact quantum group $G$, we have $C(\partial_\FH D(G))^G=\C1$ if and only if the quantum dimension function of $G$ is weakly amenable.
\end{corollary}

\bp
The ``if'' direction follows immediately from Theorem~\ref{thm:FH-weak-amenable}, as $C(\partial_\Pi\hat G)^G=\C1$.

Conversely, assume $C(\partial_\FH D(G))^G=\C1$. Take any $D(G)$-equivariant ucp map
$$
\phi\colon \RR(\ell^\infty(\hat G))\to C(\partial_\FH D(G)).
$$
Then it maps $\ell^\infty(\Irr(G))=\RR(\ell^\infty(\hat G))^G$ into scalars, hence it defines an invariant mean on $\ell^\infty(\Irr(G))$.
\ep

\begin{corollary}\label{cor:FH-trivial}
For any compact quantum group $G$, the Furstenberg--Hamana boundary of $D(G)$ is trivial if and only if $G$ is coamenable and of Kac type.
\end{corollary}

\bp
If $G$ is coamenable and of Kac type, then $\partial_\FH D(G)$ is trivial by Theorem~\ref{thm:FH-weak-amenable}, as then the quantum dimension function of $G$  is amenable and therefore $\partial_\Pi\hat G$ is trivial.

Conversely, assume $\partial_\FH D(G)$ is trivial. Then, by the previous corollary, the quantum dimension function is weakly amenable. By Theorem~\ref{thm:FH-weak-amenable} the triviality of $\partial_\FH D(G)$ is then equivalent to triviality of $\partial_\Pi\hat G$. But this means that the  quantum dimension function of $G$ is amenable, that is, $G$ is coamenable and of Kac type.
\ep

\begin{remark}
The ``if'' direction can be proved in a more elementary way as follows. We have to show that $\C$ is a $D(G)$-injective C$^*$-algebra. For this it suffices to show that the $D(G)$-injective C$^*$-algebra $\RR(\ell^\infty(\hat G))$ has a $D(G)$-invariant state. To construct such a state, we can start with any right $\hat G$-invariant mean on $\ell^\infty(\hat G)$, compose it with the left $G$-invariant conditional expectation $\ell^\infty(\hat G)\to\ell^\infty(\Irr(G))$ and then restrict the composition to $\RR(\ell^\infty(\hat G))$, see the proof of~\cite{MR2113893}*{Lemma~7.1}.

For the ``only if'' direction we can also argue as follows. If $\partial_\FH D(G)$ is trivial, then $G$ has weakly amenable quantum dimension function and $H^\infty(\hat G)=\C1$. Corollary~\ref{cor:h-infty-projection} implies then that $G$ is coamenable, while \cite{MR1916370}*{Corollary~3.9} shows that $G$ must be of Kac type. But if $\Irr(G)$ is uncountable, this argument still relies on Theorem~\ref{thm:NY-extension}. It would be interesting to find a more direct proof of such a basic property.
\hfill$\diamondsuit$
\end{remark}

For the time being it seems the only examples of compact quantum groups with weakly amenable quantum dimension functions are the quantum groups that are monoidally equivalent to coamenable compact quantum groups. For such quantum groups the noncommutative space $\partial_\FH D(G)=\partial_\Pi\hat G$ has a more explicit description.

Let us start with the coamenable case. Recall from~\cite{MR2210362} (see also~\cite{MR3556413}*{Section~2.3}) that every compact quantum group $G$ has the largest closed quantum subgroup $K$ of Kac type, namely, $\C[K]$ is the quotient of $\C[G]$ by the ideal generated by the elements $a-S^2(a)$, $a\in\C[G]$.

\begin{proposition}\label{prop:FH-coamenable}
Assume $G$ is a coamenable compact quantum group and $K$ is its maximal quantum subgroup of Kac type. Then we have the following isomorphisms of Yetter--Drinfeld $G$-C$^*$-algebras:
$$
C(\partial_\FH D(G))\cong C(\partial_\Pi\hat G)\cong C(G/K).
$$
\end{proposition}

\bp
The first isomorphism holds by Theorem~\ref{thm:FH-weak-amenable}. If $\Irr(G)$ is countable, then the second isomorphism is just a reformulation of~\cite{MR3556413}*{Theorem 3.1}; together with the last part of Theorem~\ref{thm:NY-extension} it recovers the description of $H^\infty(\hat G,\mu)$ for ergodic $\mu$ given by Tomatsu~\cite{MR2335776}. The general case is similar %, as already mentioned in~\cite{MR3556413}*{Section~3},
and can be dealt with as follows.

It is clear from Example~\ref{ex:aU} that the dimension function on $\Rep G$ defined by $C(G/K)$ coincides with the classical dimension function $U\mapsto\dim H_U$, so it is amenable by assumption. By the universality of $C(\partial_\Pi\hat G)$ we can therefore identify $C(\partial_\Pi\hat G)$ with a Yetter--Drinfeld $G$-C$^*$-subalgebra of $C(G/K)$. But by~\cite{MR3291643}*{Theorem~3.1} (which removes an extra assumption in~\cite{MR2335776}*{Theorem~3.18}) any such C$^*$-subalgebra has the form $C(G/H)$ for some intermediate closed quantum subgroup $K\subset H\subset G$. Since $C(G/H)$ defines the classical dimension function on $\Rep G$, $H$ must be of Kac type. Hence $H=K$ by the maximality of $K$.
\ep

\begin{remark}
Proposition~\ref{prop:FH-coamenable} covers all cases where $\partial_\FH D(G)$ has the form $G/H$, that is, if $G$ is a compact quantum group such that $C(\partial_\FH D(G))\cong C(G/H)$ for a closed quantum subgroup~$H$ of~$G$, then~$G$ is coamenable and $H$ is its maximal quantum subgroup $K$ of Kac type. Indeed (cf.~\cite{MR3291643}*{Proposition~4.3}), if $\partial_\FH D(G)\cong G/H$, then $C(\partial_\FH D(G))^G=\C1$, so $G$ has weakly amenable quantum dimension function. Hence $C(G/H)\cong C(\partial_\Pi G)$, which implies that $H$ is of Kac type and the classical dimension function on $\Rep G$ is amenable, so that $G$ is coamenable. As $C(G/H)\cong C(\partial_\Pi\hat G)$ has no nontrivial $D(G)$-equivariant endomorphisms, we must have $H=K$.%\hfill$\diamondsuit$
\end{remark}

\begin{example}\label{ex:Lie}
Let $G$ be a compact connected semisimple Lie group and consider its $q$-de\-for\-ma\-tion~$G_q$, $q>0$. Then by~\cite{MR2335776}*{Lemma~4.10}, for $q\ne1$, the maximal quantum subgroup of~$G_q$ of Kac type is the nondeformed maximal torus $T\subset G_q$. Therefore $\partial_\FH D(G_q)=G_q/T$ (for $q\ne1$). Since $D(G_q)$ should be thought of as a quantization of the complexification~$G_\C$ of~$G$, this agrees with \citelist{\cite{MR146298}\cite{MR161942}} showing that the Furstenberg boundary of $G_\C$ is $G_\C/P=G/T$, where $P$ is a minimal parabolic subgroup of $G_\C$. Note also that $\partial_\FH D(G)$ is trivial, as $G$ is coamenable and of Kac type. \hfill$\diamondsuit$
\end{example}

If we have a compact quantum group $G$ that is monoidally equivalent to a coamenable one,~$G_0$, then, under the correspondence between braided-commutative Yetter--Drinfeld C$^*$-algebras and tensor functors~\cite{MR3291643}, the algebras $C(\partial_\Pi\hat G)$ and $C(\partial_\Pi\hat G_0)$ correspond to the same functor from $\Rep G\sim\Rep G_0$, since the defining properties (i) and (ii) in Theorem~\ref{thm:NY-extension} can be formulated at the categorical level (see also Section~\ref{ssec:cat-FH} below). Therefore by the results of \cite{MR3291643}*{Section~3.2} we get the following generalization of Proposition~\ref{prop:FH-coamenable}.

\begin{proposition} \label{prop:FH-DRVV}
Assume $G$ is a compact quantum group that is monoidally equivalent to a coamenable compact quantum group $G_0$, and $B(G,G_0)$ is a $G$-$G_0$-Galois object defining such an equivalence. Let $K_0$ be the maximal quantum subgroup of $G_0$ of Kac type. Then we have the following isomorphisms of Yetter--Drinfeld $G$-C$^*$-algebras:
$$
C(\partial_\FH D(G))\cong C(\partial_\Pi\hat G)\cong B(G,G_0)^{K_0}.
$$
\end{proposition}

When $\Irr(G)$ is countable, the second isomorphism is basically a reformulation of a result of De Rijdt and Vander Vennet~\cite{MR2664313}*{Theorem~9.3}.

We remind (see \cite{MR3291643}*{Section~3.2}) that the action of $\hat G$ on $B(G,G_0)^{K_0}$ is given by the Miyashita--Ulbrich action. Namely, for the subalgebra of regular elements $\B\subset B(G,G_0)$ with respect to the action $\alpha\colon B(G,G_0)\to C(G)\otimes B(G,G_0)$ of $G$ we have a bijective Galois map
$$
\Gamma\colon \B\otimes_\algg\B\to\C[G]\otimes_\algg\B,\quad a\otimes b\mapsto \alpha(a)(1\otimes b).
$$
Then
$$
x\rhd a=\Gamma^{-1}(x\otimes1)_1a\Gamma^{-1}(x\otimes1)_2\quad\text{for}\quad x\in\C[G],\ a\in B(G,G_0),
$$
where we again use the sumless notation $\Gamma^{-1}(x\otimes1)=\Gamma^{-1}(x\otimes1)_1\otimes\Gamma^{-1}(x\otimes1)_2$.
The action of $\hat G$ can also be defined by a multiplicative unitary type formula, see~\cite{MR2400727}*{Eq.~(5.2)}.

\begin{corollary}
Under the assumptions of Proposition~\ref{prop:FH-DRVV}, all $D(G)$-boundaries up to isomorphism are the Yetter--Drinfeld $G$-C$^*$-algebras $B(G,G_0)^{H_0}$ for the intermediate closed quantum subgroups $K_0\subset H_0\subset G_0$.
\end{corollary}

\bp
By Theorem~\ref{thm:boundaries} all $D(G)$-boundaries up to isomorphism are the unital Yetter--Drinfeld $G$-C$^*$-subalgebras of $C(\partial_\FH D(G))$. Hence the corollary follows from Proposition~\ref{prop:FH-DRVV} and~\cite{MR3291643}*{Theorem~3.2}.
\ep

\begin{example}\label{ex:orthogonal}
For $N\ge 2$ and $Q\in\Mat_N(\C)$ such that $Q\bar Q=\pm1$, consider the free orthogonal quantum group $O^+_Q$. It is monoidally equivalent to $\SU_q(2)$, where $q\in[-1,1]\setminus\{0\}$ is determined~by
$$
\Tr(Q^*Q)=|q+q^{-1}|,\qquad \sign(Q\bar Q)=-\sign(q).
$$
An $O^+_Q$-$\SU_q(2)$-Galois object $B(O^+_Q,\SU_q(2))$ defining such an equivalence is given in~\cite{MR2202309}*{Theorem~5.5} (with $F_1=Q$ and $F_2=\begin{pmatrix}                                                                                                                                   0 & -q \\1 & 0\end{pmatrix}$). If $q\ne\pm1$, then $\T\subset \SU_q(2)$ is the maximal quantum subgroup of Kac type, hence
$$
C(\partial_\FH D(O^+_Q))\cong B(O^+_Q,\SU_q(2))^\T.
$$
Furthermore, by~\cite{MR1331688}*{Theorem~2.1} the only intermediate closed quantum subgroups $\T\subset H\subset \SU_q(2)$ are $H=\T$ and $H=\SU_q(2)$. Hence $B(O^+_Q,\SU_q(2))^\T$ is the only nontrivial $D(O^+_Q)$-boundary. We remark that, thanks to \cite{MR2400727}*{Theorem~6.1} and \cite{MR2355067}*{Theorem~5.8}, $B(O^+_Q,\SU_q(2))^\T$ has an alternative description as the space of ends of the quantum Cayley graph of ${\mathbb F}O_Q:=\widehat{ Q^+_Q}$. Note also that by~\cite{KKSV}*{Corollary~7.3}, if  $N\ge3$, then $B(O^+_Q,\SU_q(2))^\T$ is an ${\mathbb F}O_Q$-boundary.

If $q=\pm1$, then $O^+_Q\cong \SU_{\pm1}(2)$ is coamenable and of Kac type, hence the Furstenberg--Hamana boundary of $D(O^+_Q)$ is trivial. %\hfill$\diamondsuit$
\end{example}

\begin{example}\label{ex:q-auto}
Let $C$ be a finite dimensional C$^*$-algebra of dimension $\ge4$ and $\omega$ be a faithful state on $C$ such that for the multiplication map $m\colon C\otimes C\to C$ we have $mm^*=\delta^2\iota$ for some $\delta\ge(\dim C)^{1/2}$, where $m^*$ is the adjoint of $m$ with respect to the scalar product on $C$ defined by $\omega$. Then it follows from \cite{MR2664313}*{Theorem~4.7} that the quantum automorphism group $\QAut(C,\omega)$ of $(C,\omega)$ is monoidally equivalent to
$$
\SO_{q^2}(3)\cong \SU_q(2)/\{\pm1\}\cong \QAut(\Mat_2(\C),(q+q^{-1})^{-1}\Tr(\cdot\begin{pmatrix}                                                                                                                                   q& 0 \\0 & q^{-1}\end{pmatrix})),
$$
where $q\in(0,1]$ is determined by
$$
\delta=q+q^{-1}
$$
and the first isomorphism is from~\cite{MR1344000}*{Corollary~2.3}. Let $B(\QAut(C,\omega),\SU_q(2)/\{\pm1\})$ be the bi-Galois object as described in~\cite{MR2664313}*{Theorem~4.7}. If $q\ne1$, then it follows from \cite{MR1331688}*{Theorem~3.5} that the torus $\T/\{\pm1\}\cong\T$ is the maximal quantum subgroup of Kac type in $\SU_q(2)/\{\pm1\}$ and there are no intermediate quantum subgroups $\T\subset H\subset \SU_q(2)/\{\pm1\}$ apart from $H=\T$ and $H=\SU_q(2)/\{\pm1\}$. Therefore
$$
C(\partial_\FH D(\QAut(C,\omega)))\cong B(\QAut(C,\omega),\SU_q(2)/\{\pm1\})^\T
$$
and, moreover, $B(\QAut(C,\omega),\SU_q(2)/\{\pm1\})^\T$ is the only nontrivial $D(\QAut(C,\omega))$-boundary.

If $q=1$, then either $C$ is $\C^4$ and $\omega$ is given by the uniform probability distribution, or~$C$ is $\Mat_2(\C)$ and $\omega$ is the normalized trace. Therefore $\QAut(C,\omega)$ is either the quantum permutation group $S^+_4$ or the group $\SO(3)$. In either case we get a coamenable quantum group of Kac type, so the Furstenberg--Hamana boundary of $D(\QAut(C,\omega))$ is trivial.
%\hfill$\diamondsuit$
\end{example}

\subsection{Free unitary quantum groups}\label{ssec:free-unitary}

Take a natural number $N\ge 2$ and a matrix $F\in\operatorname{GL}_N(\C)$ such that $\Tr(F^*F)=\Tr((F^*F)^{-1})$. Recall that the compact free unitary quantum group $\AUF$ is defined as the universal unital C$^*$-algebra with generators $u_{ij}$, $1\le i,j\le N$, such that the matrices $U=(u_{ij})_{i,j}$ and $FU^cF^{-1}$ are unitary, where $U^c=(u_{ij}^*)_{i,j}$, equipped with the comultiplication
$$
\Delta(u_{ij})=\sum^N_{k=1}u_{ik}\otimes u_{kj}.
$$
The dual discrete quantum group is denoted by $\FAUF$.

The set $I=\Irr(\AUF)$ of isomorphism classes of irreducible representations of $\AUF$ is the free monoid on letters $\alpha$ and $\beta$, with $\alpha$ corresponding to $U$, $\beta$ to $\bar U$ and the unit $e$ to the trivial representation. The involution $s\mapsto\bar s$ is the anti-automorphism of the monoid defined by $\bar\alpha=\beta$ and $\bar\beta=\alpha$.

We will use the conventions of~\cite{MR2660685} in that we write $x$ instead of $U_x$ for $x\in I$ whenever convenient. The fusion rules for the representations of $\AUF$ are given by
$$
x\otimes y\cong\bigoplus_{z\in I:x=x_0z,y=\bar zy_0}x_0y_0.
$$
Therefore if the last letter of $x$ is the same as the first letter of $y$, then $U_x\otimes U_y$ is irreducible and isomorphic to $U_{xy}$.

The Woronowicz character $\rho$ of $\AUF$ is determined by the property
$$
\rho_U=(F^*F)^t\qquad (\text{the transpose of}\ \ F^*F).
$$
Hence $\dim_q U=\Tr(F^*F)\ge N$, and the equality holds if and only if $F$ is unitary. Let $q\in(0,1]$ be such that
\begin{equation*} \label{eq:q}
\Tr(F^*F)=q+q^{-1}.
\end{equation*}
Then $q=1$ if and only if $F$ is a unitary $2$-by-$2$ matrix.

Assume from now on that $F$ is not a unitary $2$-by-$2$ matrix, so that $q<1$. Consider the tree with vertex set $I$ such that different elements $x$ and $y$ of $I$ are connected by an edge if and only if one of them is obtained from the other by adding or removing one letter on the left. Denote by $\bar I$ the end compactification of~$I$. The elements of~$\bar I$ are words in $\alpha$ and $\beta$ that are either finite or infinite on the left, and the boundary $\partial I=\bar I\setminus I$ is the set of infinite words. The algebra $C(\bar I)$ of continuous functions on $\bar I$ can be identified with the algebra of functions $f\in\ell^\infty(I)$ such that
$$
|f(yx)-f(x)|\to0\ \text{as}\ x\to\infty,\ \text{uniformly in}\ y\in I.
$$

In \cite{MR2660685}, Vaes and Vander Vennet extended this construction to $\ell^\infty(\FAUF)$ as follows.
(To be more precise, they consider words infinite on the right, while in order to be consistent with our conventions, we consider words infinite on the left.)

For all $x,y\in I$, fix an isometry $V(yx,y\otimes x)\in\Mor(yx,y\otimes x)$. Define ucp maps
$$
\psi_{yx,x}\colon B(H_x)\to B(H_{yx})\ \ T\mapsto V(yx,y\otimes x)^*(1\otimes T)V(yx,y\otimes x).
$$
They do not depend on any choices. Define
\begin{equation*}\label{eq:VVVboundary}
C(\overline{\FAUF})=\{a\in\ell^\infty(\FAUF) : \|a_{yx}-\psi_{yx,x}(a_x)\|\to0\ \text{as}\ x\to\infty,\ \text{uniformly in}\ y\in I\}.
\end{equation*}
By \cite{MR2660685}*{Theorem~3.2}, this is a unital C$^*$-subalgebra of $\ell^\infty(\FAUF)$ containing $c_0(\FAUF)$. It can therefore be considered as the algebra of continuous functions on a (noncommutative) compactification of~$\FAUF$. The corresponding boundary is then defined by $C(\partial\FAUF)=C(\overline{\FAUF})/c_0(\FAUF)$. The Yetter--Drinfeld structure on $\RR(\ell^\infty(\FAUF))$ defines such a structure on $C(\overline{\FAUF})$ and $C(\partial\FAUF)$.

Take any generating finitely supported probability measure $\mu$ on $I$. Then by~\cite{MR2660685} the Poisson integral $\PP_\omega\colon C(\partial\FAUF)\to\ell^\infty(\FAUF)$ extends to an isomorphism of $\pi_\omega(C(\partial\FAUF))''$ onto $H^\infty(\FAUF,\mu)$, where $\omega$ is the state on $C(\overline{\FAUF})$ vanishing on $c_0(\FAUF)$ that is the weak$^*$ limit of the states $\phi_\mu^{*n}|_{C(\overline{\FAUF})}$. By~\cite{MR3961570} we can further identify $C(\partial\FAUF)$ with the Martin boundary of $\FAUF$, but we will need only the classical precursor of this result, that the Martin boundary of the random walk defined by $P_\mu|_{\ell^\infty(I)}$ is $\partial I$.

By~\cite{MR2660685}*{Proposition~4.1}, the state $\omega$ is faithful on $C(\partial{\FAUF})$, hence $\PP_\omega$ defines an isomorphism of $C(\partial\FAUF)$ onto a dense (in the ultraweak operator topology) Yetter--Drinfeld $\AUF$-C$^*$-subalgebra of $H^\infty(\FAUF,\mu)$. Note that this subalgebra is strictly smaller than $\RR(H^\infty(\FAUF,\mu))$, since the latter algebra contains a copy of $L^\infty(\partial I,\nu)$, where $\nu$ is the measure on $\partial I$ defined by~$\omega|_{C(\partial I)}$. It is probably unrealistic to have an explicit description of the Furstenberg--Hamana boundary in this case, but we at least have the following result.

\begin{theorem}\label{thm:free-unitary}
For any free unitary quantum group $\AUF$, with $F$ not a unitary $2$-by-$2$ matrix, $C(\partial\FAUF)$ is a $D(\AUF)$-boundary.
\end{theorem}

By Remark~\ref{rem:stationary}, in order to prove this theorem it suffices to show that $\omega$ is the unique $\AUF$-invariant $\phi_\mu$-stationary state on $C(\partial{\FAUF})$, cf.~\cite{KKSV}*{Section~7.2}.

As $C(\overline{\FAUF})\cap\ell^\infty(I)=C(\bar I)$ by construction, we have $C(\partial\FAUF)^{\AUF}=C(\partial I)$. It follows that every $\AUF$-invariant state on $C(\partial\FAUF)$ is determined by its restriction to $C(\partial I)$. This implies that $\phi_\mu$-stationarity should be possible to formulate entirely in terms of such restrictions. Such a reformulation is not completely straightforward though, since the operator $(\iota\otimes\phi_\mu)\Dhat$ does not leave $C(\bar I)$ invariant unless $\AUF$ is of Kac type, that is, unless~$F$ is unitary. In order to deal with this, define the states
$$
\psi_x=(\dim_q U_x)^{-1}\Tr(\cdot\,\pi_x(\rho))\quad\text{on}\quad B(H_x),
$$
and put $\psi_\mu=\sum_{x\in I}\mu(x)\psi_x$. The states $\psi_x$ are invariant with respect to the right action $T\mapsto U_x(T\otimes1)U^*_x$ of $\AUF$ on $B(H_x)$, and as a consequence the Markov operators $Q_x=(\iota\otimes\psi_x)\Dhat$ are right $\AUF$-equivariant. By restriction they define operators on~$C(\bar I)$.

Given a probability measure $\lambda$ on $\bar I$, we denote by $\lambda*\delta_x$ the measure corresponding to the state $\int Q_x(\cdot)d\lambda$ on $C(\bar I)$. We let $Q_\mu=\sum_x\mu(x)Q_x$ and $\lambda*\mu=\sum_x\mu(x)\lambda*\delta_x$.

\begin{lemma}
If $\psi$ is a $\AUF$-invariant $\phi_\mu$-stationary state on $C(\overline{\FAUF})$ and $\lambda$ is the measure on~$\bar I$ defined by~$\psi|_{C(\bar I)}$, then $\lambda$ is $\mu$-stationary, that is, $\lambda*\mu=\lambda$.
\end{lemma}

\bp
Put $Q_{\phi_\mu}=(\iota\otimes\phi_\mu)\Dhat$. As $\psi=\psi Q_{\phi_\mu}$ by assumption, it suffices to show that $\psi Q_{\phi_\mu}=\psi Q_\mu$ on $C(\bar I)$. We can approximate $\psi$ in the weak$^*$ topology by the restrictions of normal left $\AUF$-invariant states on $\ell^\infty(\FAUF)$. It follows that it is enough to show that $\phi_xQ_{\phi_\mu}=\phi_x Q_\mu$ on $\ell^\infty(I)$ for all $x\in I$. We have
$$
\phi_xQ_{\phi_\mu}=\phi_\mu P_x,\qquad \phi_x Q_\mu=\psi_\mu P_x.
$$
This gives us what we need, since $P_x$ leaves $\ell^\infty(I)$ invariant and $\phi_\mu=\psi_\mu$ on $\ell^\infty(I)$, cf.~\cite{MR2200270}*{Proposition~3.1}.
\ep

As a consequence, in order to prove Theorem~\ref{thm:free-unitary} it suffices to establish the following.

\begin{proposition}\label{prop:free-stationary2}
For any generating finitely supported probability measure $\mu$ on $I$, the measure $\nu=\lim_n\mu^{*n}$ is the unique $\mu$-stationary probability measure on $\partial I$.
\end{proposition}

The proof follows the familiar strategy for random walks on free and, more generally, hyperbolic groups, see~\cite{MR1815698}*{Theorem~2.4}. The key point is the following result.

\begin{lemma}\label{lem:convergence}
Assume $\lambda$ is a probability measure on $\bar I$ and $\{x_n\}_n$ is a sequence in $I$ converging to $x\in\partial I$. Then $\lambda*\delta_{x_n}\to\delta_x$ in the weak$^*$ topology on $C(\bar I)^*$.
\end{lemma}

\bp
For every $y\in I$, denote by $\Delta_y$ the set of all words in $I$ of the form $uy$. The closure $\bar\Delta_y$ of $\Delta_y$ in $\bar I$ consists of all words (finite and infinite) of the form $uy$. For $z\in\bar I$ of length $|z|\ge N$, denote by $[z]_N$ the word consisting of the last $N$ letters of $z$.
%Using this notation, the closure of $\Delta_y$ in $\bar I$ consists of all points $z$ such that $[z]_{N}=y$, where $N$ is the length of $y$.

The clopen sets $\bar\Delta_{[x]_N}$ form a neighbourhood base at $x$. Hence, in order to prove the lemma, it suffices to show that for every $N\ge0$ we have $(\lambda*\delta_{x_n})(\bar\Delta_{[x]_N})\to1$ as $n\to\infty$. We claim that if $y\in I$ is such that $[y]_{N+k}=[x]_{N+k}$ for some $N,k\ge0$, then
\begin{equation}\label{eq:measure-estimate}
(\lambda*\delta_{y})(\bar\Delta_{[x]_N})\ge 1-\frac{q^{2(k+1)}}{1-q^2},
\end{equation}
which obviously implies the required convergence.

Since the finitely supported probability measures on $I$ are weakly$^*$ dense in the probability measures on $\bar I$, it suffices to show~\eqref{eq:measure-estimate} for the measures $\lambda=\delta_z$, $z\in I$. We have
$$
\delta_z*\delta_y=\sum_{u\in I:z=z_0u,y=\bar uy_0}\frac{d(z_0y_0)}{d(z)d(y)}\delta_{z_0y_0},
$$
where $d$ denotes the quantum dimension. In the above sum, if $|y_0|\ge N$, then $z_0y_0\in \Delta_{[x]_N}$. Therefore
$$
(\delta_z*\delta_y)(\Delta_{[x]_N})\ge 1-\sum_{\genfrac{}{}{0pt}{2}{u\in I:z=z_0u,y=\bar uy_0}{|y_0|<N}}\frac{d(z_0y_0)}{d(z)d(y)}.
$$
By \cite{MR2660685}*{Eq.~(5)}, if $z=z_0u$ and $y=\bar u y_0$, we have
$$
d(z)\ge q^{-|u|}d(z_0),\qquad d(y)\ge q^{-|u|}d(y_0).
$$
We also have $d(z_0y_0)\le d(z_0)d(y_0)$. Hence, as $|y|\ge N+k$ by assumption,
$$
(\delta_z*\delta_y)(\Delta_{[x]_N})\ge 1-\sum_{\genfrac{}{}{0pt}{2}{u\in I:z=z_0u,y=\bar uy_0}{|y_0|<N}}q^{2|u|}
\ge 1-\sum^{N-1}_{i=0}q^{2(|y|-i)}\ge 1-\frac{q^{2(k+1)}}{1-q^2},
$$
proving~\eqref{eq:measure-estimate} for $\lambda=\delta_z$.
\ep

\bp[Proof of Proposition~\ref{prop:free-stationary2}]
The remaining proof is essentially identical to that of~\cite{MR1815698}*{Theorem~2.4}.
Consider the random walk on $I$ defined by $P_\mu|_{\ell^\infty(I)}$. Let $\Omega$ be the corresponding path space and $\PPP_e$ the Markov measure on $\Omega$ defined by the initial distribution~$\delta_e$, so that the push-forward of~$\PPP_e$ under the $n$-th projection $\Omega\to I$, $\underline{x}\to x_n$, is $\mu^{*n}$. The Martin boundary of the random walk is~$\partial I$, implying that $\PPP_e$-a.e. path $\underline{x}$ converges to a point $x_\infty\in\partial I$ and $\nu$ is the push-forward of~$\PPP_e$ under the map $\underline{x}\mapsto x_\infty$.

Assume that $\lambda$ is a $\mu$-stationary probability measure on $\bar I$. Then, for all $f\in C(\bar I)$ and $n\ge1$, we have
$$
\int_{\bar I} f\,d\lambda=\int_{\bar I} f\,d(\lambda*\mu^{*n})=\int_\Omega d\PPP_e(\underline{x})\int_{\bar I} f\, d(\lambda*\delta_{x_n}).
$$
By Lemma~\ref{lem:convergence} and the dominated convergence theorem, the last expression converges, as $n\to\infty$, to
$$
\int_\Omega f(x_\infty)d\PPP_e(\underline{x})=\int_{\partial I}f\,d\nu.
$$
Hence $\lambda=\nu$.
\ep

\section{Categorical perspective}\label{sec:cat}

For C$^*$-categories we follow the conventions of~\cite{neshveyev-tuset-book}. In particular, all such categories are assumed to be small and closed under subobjects and finite direct sums. The C$^*$-tensor categories are assumed to be strict.

\subsection{Categorification of equivariant maps}

The results of the previous sections can be formulated at the level of the representation categories $\Rep G$ and extended to C$^*$-tensor categories.
First we need to understand what the categorical analogues of ucp and completely isometric maps are.

\smallskip

Let $\CC$ be a rigid C$^*$-tensor category with simple unit. As in Section~\ref{ssec:weak-amenable}, the morphism spaces in $\CC$ will be usually denoted by $\CC(U,V)$, and we will write $\CC(U)$ for $\CC(U,U)$. Let $\Irr(\CC)$ be the set of isomorphism classes of simple objects in~$\CC$, and for every $s\in\Irr(\CC)$ fix a representative~$U_s$. The class of the unit object $\un$ is denoted by $e\in\Irr(\CC)$.

Recall that a right $\CC$-module category is a C$^*$-category $\D$ together with a unitary tensor functor from $\CC^{\mathrm{op}}$ into the category of unitary endofunctors of $\D$. Equivalently, we have a unitary bifunctor $\otimes\colon\D\times\CC\to\D$ and unitary isomorphisms $X\otimes\un\cong X$ and $(X\otimes U)\otimes V\cong X\otimes (U\otimes V)$ satisfying standard axioms, see, e.g., \cite{MR3121622}*{Section~2.3}. Unless explicitly stated otherwise, we will consider strict module categories, that is, the last two isomorphisms are assumed to be the identity morphisms.

We will be interested in pairs $(\D,X)$ consisting of a nonzero right $\CC$-module category $\D$ and an object $X\in\D$ that is generating in the sense that every object $Y$ in $\D$ is a subobject of $X\otimes U$ for some $U\in\CC$. In this case we will say, by slightly abusing the terminology, that the pair $(\D,X)$ is a \emph{singly generated right $\CC$-module category}.

For any compact quantum group $G$, there is a one-to-one correspondence between the isomorphism classes of unital $G$-C$^*$-algebras and the equivalence classes of singly generated $(\Rep G)$-module categories~\citelist{\cite{MR3121622}\cite{MR3426224}}. Specifically, the category $\D_A$ associated with a unital $G$-C$^*$-algebra $A$ is the category of finitely generated $G$-equivariant C$^*$-Hilbert $A$-modules and the generating object is $X_A=A$. Equivalently, and this is the picture we are going to use, $\D_A$ can be defined in the same way as the category $\CC_A$ we considered in Section~\ref{ssec:weak-amenable} for Yetter--Drinfeld algebras, but now the morphisms $\iota\otimes T$ are defined only for the morphisms $T$ in $\Rep G$, so $\D_A$ is only a right $(\Rep G)$-module category rather than a C$^*$-tensor category. In this picture the generating object is the unit $\un\in\Rep G\subset\D_A$.

Assume~$A_1$ and $A_2$ are two unital $G$-C$^*$-algebras. Consider the corresponding subalgebras $\A_i\subset A_i$ of regular elements and the $(\Rep G)$-module categories $\D_i=\D_{A_i}$. Then every $G$-equivariant linear map $\phi\colon\A_1\to\A_2$ defines linear maps
\begin{equation}\label{eq:multiplier}
\Theta_\phi\colon \D_1(U,V)\to\D_2(U,V)\quad\text{by}\quad\Theta_\phi(T)=(\phi\otimes\iota)(T)
\end{equation}
for $T\in\D_1(U,V)\subset\A_1\otimes B(H_U,H_V)$. As we will see shortly, the following gives an axiomatization of such maps.

\begin{definition} \label{def:cb-transformation}
Given a rigid C$^*$-tensor category $\CC$ with simple unit, a \emph{$\CC$-linear transformation} $\Theta\colon(\D_1,X_1)\to(\D_2,X_2)$ between two singly generated right $\CC$-module categories is a collection of linear maps
$$
\Theta_{U,V}\colon \D_1(X_1\otimes U,X_1\otimes V)\to\D_2(X_2\otimes U,X_2\otimes V)\qquad (U,V\in\CC)
$$
such that
\begin{enumerate}
  \item $\Theta_{U,V}((\iota\otimes S_1)T(\iota\otimes S_2))=(\iota\otimes S_1)\Theta_{Y,Z}(T)(\iota\otimes S_2)$ for all $U,V,Y,Z\in \CC$, $S_1\in\CC(Z,V)$, $S_2\in\CC(U,Y)$ and $T\in\D_1(X_1\otimes Y,X_1\otimes Z)$;
  \item $\Theta_{U\otimes Y,V\otimes Y}(T\otimes\iota_Y)=\Theta_{U,V}(T)\otimes\iota_Y$ for all $U,V,Y\in\CC$ and $T\in\D_1(X_1\otimes U,X_1\otimes V)$.
\end{enumerate}
Such a transformation is called a \emph{cb transformation}, if the maps $\Theta_{U,V}$ are bounded and
$$
\|\Theta\|_{\mathrm{cb}}:=\sup_{U,V\in\CC}\|\Theta_{U,V}\|<\infty.
$$
It is called a \emph{cp} (resp., \emph{ucp}, \emph{completely isometric}) \emph{transformation}, if the maps $\Theta_U:=\Theta_{U,U}$ are positive (resp., unital positive, isometric).
\end{definition}

We will usually write $\Theta(T)$ instead of $\Theta_{U,V}(T)$.

\smallskip

Particular cases of $\CC$-linear transformations have appeared, under different names, in~\cite{MR3679617} and~\cite{MR3406647}, see Examples~\ref{ex:C-linear} and~\ref{ex:C-linear2} below. In full generality the $\CC$-linear transformations have been introduced in~\cite{MR3687214}*{Definition~34} under the name of \emph{multipliers}, which in our opinion should rather be reserved to some special cases.

\begin{remark}\label{rem:transform}\mbox{\ }

\noindent
(1) For nonstrict module categories the definition is basically the same, but we have to use associativity morphisms in $\D_1$ and $\D_2$ to make sense of condition (ii).

\smallskip

\noindent
(2) Condition (i) means that the maps $\Theta_{U,V}$ are natural in $U$ and $V$ in the sense that they define a natural transformation between the bifunctors $\CC\times\CC\to\mathrm{Set}$, $(U,V)\mapsto\D_i(X_i\otimes U,X_i\otimes V)$ ($i=1,2$).

\smallskip

\noindent
(3) Since $\D_i(X_i\otimes U,X_i\otimes V)\otimes\Mat_n(\C)\cong \D_i(X_i\otimes U^n,X_i\otimes V^n)$, each of the spaces $\D_i(X_i\otimes U,X_i\otimes V)$ has an operator space structure. It is then clear that if $\Theta$ is a cb transformation, then the maps $\Theta_{U,V}$ are completely bounded and $\|\Theta\|_{\mathrm{cb}}=\sup_{U,V\in\CC}\|\Theta_{U,V}\|_{\mathrm{cb}}$. Similarly, if $\Theta$ is a completely isometric transformation, then the maps $\Theta_{U,V}$ are completely isometric.

Note that the same operator space structure on $\D_i(X_i\otimes U,X_i\otimes V)$ is defined by identifying $\D_i(X_i\otimes U,X_i\otimes V)$ with a corner of the C$^*$-algebra $\D_i(X_i\otimes(U\oplus V))$. If $U=V$, then this operator space structure also coincides with the one defined by the C$^*$-algebra structure on $\D_i(X_i\otimes U)$. If $\Theta$ is  a cp transformation, then the maps $\Theta_U$ are cp.

\smallskip

\noindent
(4) Any $\CC$-linear transformation $\Theta$ is determined by the maps $\Theta_{U}\colon \D_1(X_1\otimes U)\to\D_2(X_2\otimes U)$, $U\in\CC$. If we use Frobenius reciprocity and decompose into simple objects in $\CC$, then we see that $\Theta$ is also determined by the maps $\Theta_{\un,U_s}\colon \D_1(X_1,X_1\otimes U_s)\to\D_2(X_2,X_2\otimes U_s)$, $s\in\Irr(\CC)$. Furthermore, it is not difficult to show that every collection of linear maps $\D_1(X_1,X_1\otimes U_s)\to\D_2(X_2,X_2\otimes U_s)$, $s\in\Irr(\CC)$, arises this way, cf.~\cite{MR3406647}*{Proposition~3.6} and~\cite{MR3687214}*{Proposition~7}; for the representation categories this will also follow from Proposition~\ref{prop:multipliers}.

\smallskip

\noindent
(5) If a unital $\CC$-linear transformation $\Theta$ respects composition of morphisms and involution, then it gives rise to a unitary $\CC$-module functor $\D_1\to\D_2$. Namely, by replacing $\D_1$ by an equivalent category we may assume that the objects $X_1\otimes U$ ($U\in\CC$) are all different. If we denote by $\tilde\D_1$ the full subcategory of $\D_1$ consisting of such objects, then $\Theta$ defines a strict unitary $\CC$-module functor $\tilde\D_1\to\D_2$ such that $X_1\otimes U\mapsto X_2\otimes U$. This functor can then be extended to  a unitary $\CC$-module functor $\D_1\to\D_2$, since by assumption every object of $\D_1$ is a subobject of~$\tilde\D_1$. Conversely, any unitary $\CC$-module functor $\D_1\to\D_2$ mapping $X_1$ into $X_2$ defines a unital $\CC$-linear transformation that respects composition of morphisms and involution.
%\hfill$\diamondsuit$
\end{remark}

\begin{proposition}\label{prop:multipliers}
Assume $G$ is a compact quantum group, $A_1$ and $A_2$ are unital $G$-C$^*$-algebras. Consider the corresponding subalgebras $\A_i\subset A_i$ of regular elements and the right $(\Rep G)$-module categories $\D_i=\D_{A_i}$. Then the map $\phi\mapsto\Theta_\phi$ defined by~\eqref{eq:multiplier} gives a one-to-one correspondence between the $G$-equivariant linear maps $\phi\colon\A_1\to\A_2$ and the $(\Rep G)$-linear transformations $(\D_1,\un)\to (\D_2,\un)$. Furthermore, $\phi$ is cb (resp., cp, ucp, completely isometric) if and only if $\Theta_\phi$ is, and we have
\begin{equation}\label{eq:cb-norm}
\|\phi\|_\cb=\|\Theta_\phi\|_\cb.
\end{equation}
\end{proposition}

To be precise, we say that a $G$-equivariant linear map $\phi\colon\A_1\to\A_2$ is cp, if it has a (necessarily unique) extension to a cp map $A_1\to A_2$. It can be shown that this is equivalent to requiring $(\phi\otimes\iota)(T^*T)$ to be positive in $A_2\otimes\Mat_n(\C)$ for all $T\in\A_1\otimes\Mat_n(\C)$, but we won't use this.

\bp[Proof of Proposition~\ref{prop:multipliers}]
We have linear isomorphisms
$$
\pi_i\colon\bigoplus_{s\in\Irr(G)}\bar H_s\otimes\D_i(\un,U_s)\to\A_i,\quad\bar\xi\otimes T\mapsto (\iota\otimes\bar\xi)(T),
$$
for $\xi\in H_s$ and $T\in\D_i(\un,U_s)\subset\A_i\otimes B(\C,H_s)=\A_i\otimes H_s$, see~\cite{MR3426224}*{Section~2}. Then
$$
\phi\pi_1(\bar\xi\otimes T)=\pi_2(\bar\xi\otimes\Theta_\phi(T)).
$$
This shows that the map $\phi\mapsto\Theta_\phi$ is injective. This also implies that given any $(\Rep G)$-linear transformation $\Theta\colon(\D_1,\un)\to (\D_2,\un)$, we can define a $G$-equivariant linear map $\phi\colon\A_1\to\A_2$ such that $\Theta_\phi=\Theta$ on $\D_1(\un,U_s)$ for all $s\in\Irr(G)$. But then $\Theta_\phi=\Theta$ by Remark~\ref{rem:transform}(4). Therefore the map $\phi\mapsto\Theta_\phi$ is indeed a bijection.

\smallskip

Next, it is clear by definition that if $\phi$ is cb, then $\Theta_\phi$ is also cb and $\|\Theta_\phi\|_\cb\le\|\phi\|_\cb$. Assume now that $\Theta_\phi$ is cb. Let us show first that $\phi$ is bounded and $\|\phi\|\le\|\Theta_\phi\|_\cb$. Assume the $G$-action on $A_i$ is given by $\alpha_i\colon A_i\to C(G)\otimes A_i$.
Consider the right regular representation $V\in M(K(L^2(G))\otimes C(G))$ of $G$.
%Namely, $V=(\hat J\otimes\hat J)W^*_{21}(\hat J\otimes\hat J)$.
It has the property that
$$
V_{31}^*(\alpha_i\otimes\iota)(T)V_{31}=T\quad\text{for all}\quad T\in\alpha_i(A_i)_{21}\in A_i\otimes C(G),
$$
where we view the last factor $C(G)$ as a subalgebra of $B(L^2(G))$. We can find a net of finite rank $G$-invariant projections $p_j\in B(L^2(G))$ converging strongly to~$1$. Let $V_j\in\Rep G$ be the restriction of $V$ to $p_jL^2(G)$. Then the above identity implies that
$$
T_{ij}(a):=(1\otimes p_j)\alpha_i(a)_{21}(1\otimes p_j)\in\D_i(V_j)\quad\text{for all}\quad a\in A_i.
$$
By the $G$-equivariance of $\phi$ we also have
\begin{equation}\label{eq:equiv}
\Theta_\phi(T_{1j}(a))=T_{2j}(\phi(a))\quad\text{for all}\quad a\in\A_1,
\end{equation}
and therefore $\|T_{2j}(\phi(a))\|\le\|\Theta_\phi\|_\cb\|T_{1j}(a)\|$. Passing to the limit we get
$\|\phi(a)\|\le\|\Theta_\phi\|_\cb\|a\|$ for all $a\in\A_1$, as required.

Repeating the above argument with $\A_i$ replaced by $\A_i\otimes\Mat_n(\C)$ and $V$ by $V\oplus\dots\oplus V$, we get that $\|\phi\|_\cb\le \|\Theta_\phi\|_\cb$, finishing the proof of~\eqref{eq:cb-norm}.

\smallskip

The same arguments as above show that $\phi$ is completely isometric if and only if $\Theta_\phi$ is. It is also obvious that $\phi$ is unital if and only if $\Theta_\phi$ is.

\smallskip

It remains to deal with complete positivity. It is again clear by definition that if $\phi$ is cp, then $\Theta_\phi$ is also cp. Assume now that $\Theta_\phi$ is cp. Then from the cb case we can already conclude that $\phi$ is cb and $\|\phi\|_\cb=\|\Theta_\phi\|_\cb=\|\phi(1)\|$. Hence $\phi$ extends to a cb map $A_1\to A_2$, which we continue to denote by $\phi$. Let us show that $\phi$ is positive. This can again be deduced from~\eqref{eq:equiv}, but it is also possible to give a more direct algebraic proof as follows.

Take $a\in \A_1$. Then $a$ is the sum of finitely many elements $(\iota\otimes\bar\xi_k)(T_k)$, with $\xi_k\in H_{s_k}$, $T_k\in \D_1(\un,U_{s_k})$. Put $U=\bigoplus_k U_{s_k}$. Then we can define $\xi\in H_U$ and $T\in\D_1(\un,U)$ such that
$a=(\iota\otimes\bar\xi)(T)$. It follows that $aa^*=(\iota\otimes\omega_{\xi,\xi})(TT^*)$,
where $TT^*\in\D_1(U)$ and $\omega_{\xi,\xi}=(\cdot\,\xi,\xi)$ is a positive linear functional on $B(H_U)$. Hence
$$
\phi(aa^*)=(\iota\otimes\omega_{\xi,\xi})\Theta_\phi(TT^*)\ge0,
$$
so that $\phi$ is positive. Passing to the matrix algebras over $\A_i$, a similar argument shows that $\phi$ is cp.
\ep

We are, however, more interested in braided-commutative Yetter--Drinfeld $G$-C$^*$-algebras $A$. In this case, as we discussed in Section~\ref{ssec:weak-amenable}, the $\hat G$-action defines a tensor structure on $\D_A$, so that we get a C$^*$-tensor category $\CC_A$ containing $\Rep G$.

Assume again that $\CC$ is a rigid C$^*$-tensor category with simple unit. If we are given a C$^*$-tensor category $\B$ and a dominant unitary tensor functor $\CC\to\B$, then we say that $\B$ is a \emph{$\CC$-tensor category}. By replacing~$\B$ by an equivalent category we can and will tacitly assume that $\B$ is strict and $\CC$ is simply a C$^*$-tensor subcategory of $\B$.

\begin{definition} \label{def:cb-transformation2}
A \emph{$\CC$-linear transformation} $\Theta\colon\CC_1\to\CC_2$ between two $\CC$-tensor categories is a collection of linear maps
$$
\Theta_{U,V}\colon \CC_1(U,V)\to\CC_2(U,V)\qquad (U,V\in\CC)
$$
such that
\begin{enumerate}
  \item $\Theta(S_1TS_2)=S_1\Theta(T)S_2$ for all $U,V,Y,Z\in \CC$, $S_1\in\CC(Z,V)$, $S_2\in\CC(U,Y)$ and $T\in\CC_1(Y,Z)$;
  \item $\Theta(\iota_Z\otimes T\otimes\iota_Y)=\iota_Z\otimes\Theta(T)\otimes\iota_Y$ for all $U,V,Y,Z\in\CC$ and $T\in\CC_1(U,V)$.
\end{enumerate}
\end{definition}

In other words, a $\CC$-linear transformation $\Theta\colon\CC_1\to\CC_2$ is a collection of maps $\Theta_{U,V}\colon \CC_1(U,V)\to\CC_2(U,V)$ that defines $\CC$-linear transformations $(\CC_1,\un)\to(\CC_2,\un)$ if we consider $\CC_i$ separately as left and right $\CC$-module categories. Equivalently, it is a $\CC$-linear transformation $(\CC_1,\un)\to(\CC_2,\un)$ of right $\CC$-module categories satisfying the extra condition
$$
\Theta(\iota_Y\otimes T)=\iota_Y\otimes\Theta(T)
$$
for all $U,V,Y\in\CC$ and $T\in\CC_1(U,V)$.

\begin{proposition}\label{prop:multipliers2}
Assume $G$ is a compact quantum group, $A_1$ and $A_2$ are unital braided-commutative Yetter--Drinfeld $G$-C$^*$-algebras. Consider the corresponding subalgebras $\A_i\subset A_i$ of regular elements and the $(\Rep G)$-tensor categories $\CC_i=\CC_{A_i}$. Then a $G$-equivariant linear map $\phi\colon\A_1\to\A_2$ is $\hat G$-equivariant if and only if the corresponding $(\Rep G)$-linear transformation $\Theta_\phi\colon (\CC_1,\un)\to(\CC_2,\un)$ of singly generated right $(\Rep G)$-module categories is a $(\Rep G)$-linear transformation $\CC_1\to\CC_2$ of $(\Rep G)$-tensor categories.
\end{proposition}

\bp
This follows immediately from the definition of the tensor structure on $\CC_i$ in terms of the $\hat G$-action, see~\eqref{eq:tensor-structure}.
\ep

\begin{example}\label{ex:C-linear}
Given a $\CC$-tensor category $\B$, a ucp transformation $\hat\CC\to\B$ (see the next subsection for the definition of $\hat\CC$) is the same thing as a right invariant mean on the functor $\CC\to\B$ in the sense of~\cite{MR3679617}*{Definition~6.1}.%\hfill$\diamondsuit$
\end{example}

\begin{example}[cf. \cite{MR3687214}*{Example~15}]\label{ex:C-linear2}
We can view $\CC$ as a $\CC$-bimodule category, or equivalently, as a right $(\CC^{\mathrm{op}}\boxtimes\CC)$-module category. Then a $(\CC^{\mathrm{op}}\boxtimes\CC)$-linear transformation $(\CC,\un)\to(\CC,\un)$ is the same thing as a multiplier on $\CC$ in the sense of \cite{MR3406647}*{Definition~3.4} (see also~\cite{MR3406647}*{Proposition~3.6}). As an example consider the representation category $\CC=\Rep G$ of a compact quantum group $G$. Then the C$^*$-algebra corresponding to the bimodule category $\Rep G$ and the generating object $\un$ is $C(G)$ equipped with its usual left and right actions of $G$ by translations. Therefore by Proposition~\ref{prop:multipliers} we get a one-to-one correspondence between the cb-multipliers (resp., cp-multipliers) on $\Rep G$ and the $G$-$G$-equivariant cb (resp., cp) maps $C(G)\to C(G)$, thus recovering~\cite{MR3406647}*{Proposition~6.1}, but also proving equality~\eqref{eq:cb-norm}.

On the other hand, if we view $\CC$ either as a right or left $\CC$-module category with generating object $\un$ or as a $\CC$-tensor category, then $\CC(\un,U_s)=0$ for all $s\ne e$ and the only $\CC$-linear transformations $\Theta$ on $\CC$ are the scalars, that is, $\Theta_{U,V}\colon\CC(U,V)\to\CC(U,V)$ is the multiplication by a scalar, the same one for all~$U$ and~$V$. %\hfill$\diamondsuit$
\end{example}

\begin{example}
Consider a compact quantum group $G$. The forgetful functor $\Rep G\to\Hilbf$ allows us to view the tensor category $\Hilbf$ of finite dimensional Hilbert spaces as a right $(\Rep G)$-module category. The $G$-C$^*$-algebra corresponding to $(\Hilbf,\C)$ is~$C(G)$ equipped with the action of $G$ by left translations. The $G$-equivariant linear maps $\C[G]\to\C[G]$ have the form $m_\phi=(\iota\otimes\phi)\Delta$ for linear functionals $\phi$ on $\C[G]$. If $m_\phi$ is cb (resp., cp), then~$\phi$ is called a cb (resp., cp) multiplier, or a Herz--Schur multiplier, on $\hat G$. By Proposition~\ref{prop:multipliers} we therefore get a one-to-one correspondence between the cb (resp., cp) multipliers on $\hat G$ and the cb (resp., cp) $(\Rep G)$-linear transformations $(\Hilbf,\C)\to(\Hilbf,\C)$.

It is known that $\phi\in\C[G]^*$ is a cp multiplier if and only if the linear functional $\phi$ is positive, see~\cite{MR784292}*{Proposition~4.2}, \cite{MR3019431}*{Theorem~5.2}. Let us give a quick proof of this fact. If~$\phi$ is positive, then it extends to a positive linear functional on the universal completion $C_u(G)$ of~$\C[G]$. As the comultiplication $\Delta$ on $\C[G]$ extends to a $*$-homomorphism $C(G)\to C(G)\otimes C_u(G)$, it follows that $m_\phi=(\iota\otimes\phi)\Delta$ extends to a cp map on $C(G)$. Conversely, assume~$\phi$ is a cp multiplier. Take $a\in\C[G]$. As in the last part of the proof of Proposition~\ref{prop:multipliers}, we can find an element
$T\in\D_{C(G)}(\un,U)$ and a state $\omega$ on $B(H_U)$ such that $aa^*=(\iota\otimes\omega)(TT^*)$. Then $(m_\phi\otimes\iota)(TT^*)$ is a positive element of $\D_{C(G)}(U)$, so it has the form~$SS^*$ for some $S\in\D_{C(G)}(U)\subset \C[G]\otimes B(H_U)$. But then
$$
\phi(aa^*)=\eps(m_\phi(aa^*))=\omega\big((\eps\otimes\iota)(S)(\eps\otimes\iota)(S)^*\big)\ge0,
$$
so that $\phi$ is positive.

We thus get a one-to-one correspondence between the positive linear functionals on $\C[G]$ and the cp $(\Rep G)$-linear transformations $(\Hilbf,\C)\to(\Hilbf,\C)$. This recovers \cite{MR3687214}*{Proposition~15}.
\end{example}

\subsection{Furstenberg--Hamana boundaries of monoidal categories}\label{ssec:cat-FH}

We continue to assume that $\CC$ is a rigid C$^*$-tensor category with simple unit. We say that a $\CC$-tensor category $\A$ is \emph{$\CC$-injective} if, given $\CC$-tensor categories $\B$ and $\D$, a completely isometric ucp transformation $\Phi\colon \B\to \D$ and a ucp transformation $\Psi\colon \B\to \A$, there is a ucp transformation $\tilde\Psi\colon \D\to \A$ making the diagram
$$
\xymatrix{\D\ar@{-->}[rd]^{\tilde\Psi} & \\ \B\ar[u]^{\Phi}\ar[r]_{\Psi} & \A}
$$
commutative, that is, $\Psi_{U,V}=\tilde\Psi_{U,V}\Phi_{U,V}$ for all $U,V\in\CC$.

\begin{definition}
A  $\CC$-tensor category $\A$ is called a \emph{Furstenberg--Hamana boundary} of $\CC$, if $\A$ is $\CC$-injective and, for all $\CC$-tensor categories $\B$, every ucp transformation $\A\to\B$ is completely isometric.
\end{definition}

Let us introduce the following terminology. Given two $\CC$-tensor categories $\A$ and $\B$, by a \emph{$\CC$-tensor functor} $\A\to\B$ we mean a unitary tensor functor $\A\to\B$ such that its composition with the functor $\CC\to\A$ is naturally unitarily monoidally isomorphic to the functor $\CC\to\B$. If $\CC$ is a subcategory of $\A$ and $\B$, as we usually assume, then a $\CC$-tensor functor can be assumed to be identical on~$\CC$, see the discussion in~\cite{MR3291643}*{Section~2.1}. We say that $\A$ and $\B$ are equivalent as $\CC$-tensor categories, if there is a $\CC$-tensor functor $\A\to\B$ that is an equivalence of categories.

\begin{theorem}\label{thm:FH2}
For any rigid C$^*$-tensor category $\CC$ with simple unit, there is a unique up to equivalence Furstenberg--Hamana boundary of $\CC$.
\end{theorem}

We denote the Furstenberg--Hamana boundary of $\CC$ by $\partial_\FH\CC$.

\smallskip

In order to prove the theorem we need a categorical analogue of $\RR(\ell^\infty(\hat G))$. This is the $\CC$-tensor category $\hat\CC$ introduced in~\cite{MR3679617}. For $U,V\in\CC$, the morphism space $\hat\CC(U,V)$ is defined as the space $\Nat_b(\iota\otimes U,\iota\otimes V)$ of bounded natural transformations between the functors $\iota\otimes U$ and $\iota\otimes V$ of tensoring on the right on $\CC$. The tensor product of natural transformations is defined by the following rules.
Given $\nu=(\nu_X\colon X\otimes U\to X\otimes V)_{X\in\CC}\in\hat\CC(U,V)$, we have
$$
(\nu\otimes\iota_Y)_X=\nu_X\otimes\iota_Y,\qquad (\iota_Y\otimes\nu)_X=\nu_{X\otimes Y}.
$$

The following is a (partial) categorical analogue of Proposition~\ref{prop:Izumi-equivariance}.

\begin{proposition}\label{prop:Izumi-equivariance2}
For any $\CC$-tensor category $\A$, there is a one-to-one correspondence between the ucp transformations $\A\to\hat\CC$ and the states on the C$^*$-algebra $\A(\un)=\End_\A(\un)$.
\end{proposition}

We denote by $\PP_\psi$ the ucp transformation $\A\to\hat\CC$ defined by a state $\psi$ on $\A(\un)$, and call it a \emph{Poisson transformation}.

\bp
Assume $\PP\colon\A\to\hat\CC$ is a ucp transformation. Composing $\PP\colon\A(\un)\to\hat\CC(\un)\cong\ell^\infty(\Irr(\CC))$ with the evaluation at $e\in\Irr(\CC)$, we get a state $\psi$ on $\A(\un)$. Take an object $U$ in $\CC$ and choose a standard solution $(R_U,\bar R_U)$ of the conjugate equations for $U$ in $\CC$. Then we can define a state~$\psi_U$ on~$\A(U)$~by
$$
\psi_U(T)=d^\CC(U)^{-1}\psi(\bar R_U(T\otimes\iota)\bar R_U^*),
$$
where $d^\CC(U)=\|R_U\|^2=\|\bar R_U\|^2$ is the intrinsic dimension of $U$ in $\CC$. The same argument as in the proof of~\cite{MR3679617}*{Lemma~4.1} shows that, as $\psi_U$ is independent of the choice of standard solutions, $\CC(U)$ is contained in the centralizer of $\psi_U$. Since $\CC(U)$ is finite dimensional and~$\psi_U$~is faithful on it, it follows that there is a unique $\psi_U$-preserving conditional expectation $E_U\colon\A(U)\to\CC(U)$.

We claim that $\PP$ is given by
\begin{equation}\label{eq:PP(T)}
\PP(T)_X=E_{X\otimes U,X\otimes U}(\iota_X\otimes T)\colon X\otimes U\to X\otimes U
\end{equation}
for all $U,X\in\CC$, $T\in\A(U)$. In order to prove this, let us first of all observe that we have
\begin{equation*}\label{eq:psiU}
\psi_U(T)=\tr^\CC_U(\PP(T)_\un)\quad\text{for all}\quad T\in\A(U),
\end{equation*}
since $\PP(\bar R^*_U(T\otimes\iota_{\bar U})R_U)_\un=\bar R_U^*(\PP(T)_\un\otimes\iota_{\bar U})\bar R_U$, where $\tr^\CC$ is the normalized categorical trace. Since we also have
$$
\PP(T)_X=\PP(T)_{X\otimes\un}=(\iota_X\otimes\PP(T))_\un=\PP(\iota_X\otimes T)_\un,
$$
we get, for all $U,X\in\CC$, $T\in\A(U)$ and $S\in\CC(X\otimes U)$, that
$$
\tr^\CC_{X\otimes U}(S\PP(T)_X)=\tr^\CC_{X\otimes U}(S\PP(\iota_X\otimes T)_\un)
=\tr^\CC_{X\otimes U}(\PP(S(\iota_X\otimes T))_\un)=\psi_{X\otimes U}(S(\iota_X\otimes T)),
$$
which implies~\eqref{eq:PP(T)}. Therefore $\PP$ is completely determined by $\psi$.

Conversely, starting from an arbitrary state $\psi$ on $\A(\un)$ we define conditional expectations~$E_U$ as above and then define $\PP\colon\A(U)\to\hat\CC(U)$ by~\eqref{eq:PP(T)}. By identifying $\A(U,V)$ with a corner of $\A(U\oplus V)$, we also get maps $\A(U,V)\to\hat\CC(U,V)$. It is straightforward to check that $\PP$ is a ucp transformation, cf.~\cite{MR3679617}*{Lemma~4.5}.
\ep

\bp[Proof of Theorem~\ref{thm:FH2}]
Consider the convex semigroup $S$ of ucp $\CC$-linear transformations $\hat\CC\to\hat\CC$. We can think of the elements of $S$ as maps on the Banach space
$$
\ell^\infty\text{-}\bigoplus_{U,V\in\CC}\hat\CC(U,V)\cong\ell^\infty\text{-}\bigoplus_{s\in\Irr(\CC),U,V\in\CC}\CC(U_s\otimes U,U_s\otimes V),
$$
which is dual to $\ell^1\text{-}\bigoplus_{s\in\Irr(\CC),U,V\in\CC}\CC(U_s\otimes U,U_s\otimes V)^*$. It is not difficult to see that $S$ is closed in the topology of pointwise weak$^*$ convergence. Hence, by Proposition~\ref{prop:Hamana}, there is a minimal idempotent $\Theta\in S$. We define a new $\CC$-tensor category $\A$, with
$$
\A(U,V)=\Theta(\hat\CC(U,V))\quad (U,V\in\CC)
$$
and the composition of morphisms given by the Choi--Effros product $T_1\cdot T_2=\Theta(T_1T_2)$, cf.~\cite{MR3679617}*{Section~2}. The rest of the argument, showing that $\A$ is a unique up to equivalence Furstenberg--Hamana boundary of $\CC$, is similar to the proof of Theorem~\ref{thm:FH} (starting from Step 4), but now using Proposition~\ref{prop:Izumi-equivariance2} instead of Proposition~\ref{prop:Izumi-equivariance}. We skip the details.
\ep

Recall that by Corollary~\ref{cor:br-com} the Yetter--Drinfeld C$^*$-algebras $C(\partial_\FH D(G))$ are braided\--commu\-ta\-tive. Hence we have $(\Rep G)$-tensor categories associated with them.

\begin{proposition}\label{prop:FH-catFH}
For any compact quantum group $G$, the $(\Rep G)$-tensor category associated with $C(\partial_\FH D(G))$ is equivalent to $\partial_\FH(\Rep G)$.
\end{proposition}

\bp
In view of Propositions~\ref{prop:multipliers} and~\ref{prop:multipliers2}, we have a one-to-one correspondence between the $D(G)$-equivariant ucp maps from one braided-commutative Yetter--Drinfeld $G$-C$^*$-algebra into another, and the ucp transformations between the associated $(\Rep G)$-tensor categories. This implies that the $(\Rep G)$-tensor category associated with $C(\partial_\FH D(G))$ has the defining properties of $\partial_\FH(\Rep G)$, hence it is equivalent to $\partial_\FH(\Rep G)$.
\ep

This result follows also from the construction of $C(\partial_\FH D(G))$ and $\partial_\FH(\Rep G)$, since the semigroups $S$ used in both cases become literally the same, hence they have the same minimal idempotents, under the correspondence between the $D(G)$-equivariant ucp maps on $\RR(\ell^\infty(\hat G))$ and the ucp $(\Rep G)$-linear transformations on $\widehat{\Rep G}$.

\subsection{Monoidal categories with weakly amenable dimension functions}\label{ssec:weakly-amenable2}

The Markov operators $P_s$ introduced in Section~\ref{ssec:Poisson} have a categorical analogue~\cite{MR3679617}. Namely, using the notation of the previous subsection, for $\nu=(\nu_X)_X\in\hat\CC(U,V)$, we define
$$
P_s(\nu)_X=(\tr^\CC_{U_s}\otimes\iota)(\nu_{U_s\otimes X})\in\CC(X\otimes U,X\otimes V).
$$
Therefore $P_s\colon\hat\CC\to\hat\CC$ is a ucp $\CC$-linear transformation. Recall then that a probability measure~$\mu$ on~$\Irr(\CC)$ is called \emph{ergodic}, if the constants are the only $P_\mu$-harmonic elements in $\hat\CC(\un)\cong\ell^\infty(\Irr(\CC))$, where $P_\mu=\sum_{s\in\Irr(\CC)}\mu(s)P_s$.

\begin{theorem}\label{thm:NY-extension2}
Assume $\CC$ is a rigid C$^*$-tensor category with simple unit and weakly amenable intrinsic dimension function. Then there is a unique up to equivalence $\CC$-tensor category $\partial_\Pi\CC$ with simple unit such that the following conditions are satisfied:
\begin{enumerate}
\item the dimension function on $\CC$ defined by the intrinsic dimension function on $\partial_\Pi\CC$ is amenable;
\item if $\B$ is another $\CC$-tensor category with simple unit that defines the amenable dimension function on $\CC$, then there is
a $\CC$-tensor functor $\partial_\Pi\CC\to\B$, unique up to natural unitary monoidal isomorphism.
\end{enumerate}

Furthermore, if $\phi$ is the unique state on $\End_{\partial_\Pi\CC}(\un)\cong\C$, then the Poisson transformation $\PP_\phi\colon \partial_\Pi\CC\to\hat\CC$ is completely isometric and, for all $U,V\in\CC$, the image of $\Mor_{\partial_\Pi\CC}(U,V)$ under~$\PP_\phi$~is
$$
\{\nu\in {\hat\CC}(U,V)\mid P_s(\nu)=\nu\ \text{for all}\ s\in\Irr(\CC)\}.
$$
If $\Irr(\CC)$ is countable, then the last space coincides with the space of $P_\mu$-harmonic bounded natural transformations $\iota\otimes U\to\iota\otimes V$ for any ergodic probability measure~$\mu$ on~$\Irr(\CC)$.
\end{theorem}

\bp
When $\Irr(\CC)$ is countable, this is a reformulation of Theorems 4.1 and 5.1 in~\cite{MR3679617}. The result is then extended to arbitrary categories similarly to the proof of Theorem~\ref{thm:NY-extension}, by considering a net of full subcategories $(\CC_i)_i$ with countable $\Irr(\CC_i)$ and weakly amenable intrinsic dimension functions. We omit the details.
\ep

\begin{theorem} \label{thm:FH-weak-amenable2}
For any rigid C$^*$-tensor category $\CC$ with simple unit and weakly amenable intrinsic dimension function, the $\CC$-tensor categories $\partial_\FH\CC$ and $\partial_\Pi\CC$ are equivalent.
\end{theorem}

\bp
The proof is similar to that of Theorem~\ref{thm:FH-weak-amenable}. Let us show that $\partial_\Pi\CC$ is $\CC$-injective. Assume first that $\Irr(\CC)$ countable. Choose an ergodic probability measure $\mu$ on $\Irr(\CC)$, a free ultrafilter~$\omega$ on~$\N$ and define an idempotent ucp transformation $E\colon\hat\CC\to\hat\CC$ by
$$
E(\nu)_X=\lim_{n\to\omega}\frac{1}{n}\sum^n_{k=1}P^k_\mu(\nu)_X
$$
for all $\nu\in\hat\CC(U,V)$ and $X\in\CC$. Then the image of $E$ consists exactly of the $P_\mu$-harmonic bounded natural transformations, hence, by Theorem~\ref{thm:NY-extension2}, it coincides with the image of  $\PP_\phi$. As $\PP_\phi$ is completely isometric and~$\hat\CC$ is $\CC$-injective by Proposition~\ref{prop:Izumi-equivariance2}, it follows that $\partial_\Pi\CC$ is $\CC$-injective as well. The case when $\Irr(\CC)$ is uncountable is dealt with similarly to the proof of Corollary~\ref{cor:h-infty-projection}, we omit the details.

Next, assume $\Theta\colon\partial_\Pi\CC\to\B$ is a ucp transformation for some $\CC$-tensor category $\B$. Consider any ucp transformation $\Psi\colon\B\to\hat\CC$. Then $\Psi\Theta=\PP_\phi$ by Proposition~\ref{prop:Izumi-equivariance2}, as $\phi$ is the only state on $\End_{\partial_\Pi\CC}(\un)$. Thus, $\Psi\Theta$ is completely isometric, so $\Theta$ must be completely isometric as well. Therefore $\partial_\Pi\CC$ is a Furstenberg--Hamana boundary of $\CC$.
\ep

The following result is a generalization of Corollary~\ref{cor:FH-G-ergodic}.

\begin{corollary}
For any rigid C$^*$-tensor category $\CC$ with simple unit, the category $\partial_\FH\CC$ has simple unit if and only if $\CC$ has weakly amenable intrinsic dimension function.
\end{corollary}

\bp
The ``if'' direction is immediate, as $\partial_\Pi\CC$ has simple unit by definition. Conversely, assume~$\partial_\FH\CC$ has simple unit. By the $\CC$-injectivity there is a ucp transformation $\Theta\colon\hat\CC\to\partial_\FH\CC$. In the terminology of~\cite{MR3679617} this means that the functor $\CC\to\partial_\FH\CC$ is amenable. Then by~\cite{MR3679617}*{Lemma~6.4} the map~$\Theta_\un$ is an invariant mean on $\ell^\infty(\Irr(\CC))$, so that the intrinsic dimension function on $\CC$ is weakly amenable.
\ep

The following is a generalization of Corollary~\ref{cor:FH-trivial}.

\begin{corollary}
For any rigid C$^*$-tensor category $\CC$ with simple unit, the Furstenberg--Hamana boundary of $\CC$ is trivial if and only if $\CC$ is amenable.
\end{corollary}

The triviality here means that the $\CC$-tensor category $\partial_\FH\CC$ is equivalent to $\CC$.

\bp
Again, the ``if'' direction is immediate, as if $\CC$ is amenable, then $\CC$ has the defining properties of $\partial_\Pi\CC$, so that we can take $\partial_\FH\CC=\partial_\Pi\CC=\CC$. Conversely, assume $\partial_\FH\CC$ is trivial. By the previous corollary, the intrinsic dimension function on $\CC$ is weakly amenable. But then~$\partial_\Pi\CC$ is equivalent to $\partial_\FH\CC\sim\CC$, which means that the intrinsic dimension function on~$\CC$ is amenable.
\ep

The $q$-deformations provide arguably the most important examples of nonamenable C$^*$-tensor categories with weakly amenable intrinsic dimension functions. For these categories we have, thanks to Proposition~\ref{prop:FH-catFH}, the following reformulation of Examples~\ref{ex:Lie},~\ref{ex:orthogonal} and~\ref{ex:q-auto}.

\begin{example}
Assume $G$ is a compact connected semisimple Lie group and consider its $q$-deformation $G_q$, $q>0$. Let $T\subset G_q$ be the nondeformed maximal torus. Then, for $q\ne1$, we have $\partial_\FH(\Rep G_q)=\Rep T$, with the functor $\Rep G_q\to\Rep T$ being simply the forgetful functor. For $q=1$, the Furstenberg--Hamana boundary of $\Rep G$ is trivial.%\hfill$\diamondsuit$
\end{example}

\bigskip
% \bib, bibdiv, biblist are defined by the amsrefs package.

\bigskip

\end{document}